\documentclass[11pt]{amsart}
\newtheorem{theorem}{Theorem}[section]
\newtheorem{proposition}[theorem]{Proposition}
\newtheorem{lemma}[theorem]{Lemma}
\newtheorem{corollary}[theorem]{Corollary}
\theoremstyle{definition}
\newtheorem{definition}[theorem]{Definition}

\begin{document}

\title[Ancient solutions to mean curvature flow in $\mathbb{R}^3$]{Uniqueness of convex ancient solutions to mean curvature flow in $\mathbb{R}^3$}
\author{Simon Brendle and Kyeongsu Choi}
\address{Department of Mathematics, Columbia University, 2990 Broadway, New York, NY 10027, USA.}
\address{Department of Mathematics, Massachusetts Institute of Technology, 77 Massachusetts Avenue, Cambridge MA 02138, USA.}
\begin{abstract}
A well-known question of Perelman concerns the classification of noncompact ancient solutions to the Ricci flow in dimension $3$ which have positive sectional curvature and are $\kappa$-noncollapsed. In this paper, we solve the analogous problem for mean curvature flow in $\mathbb{R}^3$, and prove that the rotationally symmetric bowl soliton is the only noncompact ancient solution of mean curvature flow in $\mathbb{R}^3$ which is strictly convex and noncollapsed. 
\end{abstract}
\thanks{The first author was supported by the National Science Foundation under grant DMS-1649174 and by the Simons Foundation. The second author was supported by the National Science Foundation under grant DMS-1811267. The authors gratefully acknowledge the hospitality of T\"ubingen University, where part of this work was carried out.}
\maketitle

\section{Introduction}

This paper is concerned with the classification of ancient solutions to mean curvature flow. By definition, an ancient solution is a solution which is defined for $t \in (-\infty,T]$ for some $T$. Ancient solutions play an important role in understanding singularity formation in geometric flows. For example, Perelman's famous Canonical Neighborhood Theorem \cite{Perelman} states that, for a solution to the Ricci flow in dimension $3$, the high curvature regions are modeled on ancient solutions, which have nonnegative curvature and are $\kappa$-noncollapsed. Moreover, Perelman \cite{Perelman} proved that every noncompact ancient $\kappa$-solution in dimension $3$ has the structure of a tube with a positively curved cap attached on one side.

Results of a similar nature hold for mean curvature flow, assuming that the initial hypersurfacer is mean convex and embedded (see \cite{Brendle-Huisken},\cite{Haslhofer-Kleiner1},\cite{Haslhofer-Kleiner2},\cite{Huisken-Sinestrari1},\cite{Huisken-Sinestrari2},\cite{White1},\cite{White2}). In particular, Huisken-Sinestrari \cite{Huisken-Sinestrari1},\cite{Huisken-Sinestrari2} proved that, for any mean convex solution to mean curvature flow, the high curvature regions are almost convex. Under the stronger assumption of two-convexity, one can show that the mean curvature flow will only form neck-pinch singularities. Moreover, the flow can be continued beyond singularities by a surgery procedure similar in spirit to the one devised by Hamilton and Perelman for the Ricci flow (see \cite{Brendle-Huisken},\cite{Haslhofer-Kleiner2},\cite{Huisken-Sinestrari3}).

Our goal in this paper is to classify all convex ancient solutions to mean curvature flow in $\mathbb{R}^3$ which are $\alpha$-noncollapsed in the sense of Sheng and Wang \cite{Sheng-Wang}. Recall that a mean convex surface is called $\alpha$-noncollapsed if, for each point $x$ on the surface, there exists a ball of radius $\frac{\alpha}{H}$ in ambient space, which lies inside the surface and which touches the surface at the given point $x$. Examples of noncollapsed ancient solutions include the shrinking cylinders and the rotationally symmetric bowl soliton (cf. \cite{Altschuler-Wu}). In this paper, we show that these are the only ancient solutions which are noncompact, convex, and noncollapsed:

\begin{theorem} 
\label{main.thm}
Let $M_t$, $t \in (-\infty,0]$, be a noncompact ancient solution of mean curvature flow in $\mathbb{R}^3$ which is strictly convex and noncollapsed. Then $M_t$ agrees with the bowl soliton, up to scaling and ambient isometries. 
\end{theorem}

By combining Theorem \ref{main.thm} with known results in the literature (see \cite{Haslhofer-Kleiner1}, Theorem 1.10, or \cite{White1},\cite{White2}), we can draw the following conclusion:

\begin{corollary}
\label{singularity.formation}
Consider an arbitrary closed, embedded, mean convex surface in $\mathbb{R}^3$, and evolve it by mean curvature flow. At the first singular time, the only possible blow-up limits are shrinking round spheres; shrinking round cylinders; and the translating bowl soliton.
\end{corollary}

Let us indicate how Corollary \ref{singularity.formation} follows from Theorem \ref{main.thm}. Suppose we evolve a closed, embedded, mean convex surface in $\mathbb{R}^3$ by mean curvature flow. At the first singular time, every blow-up limit is an ancient solution which is weakly convex and noncollapsed (cf. \cite{White1}, \cite{White2}, or \cite{Haslhofer-Kleiner1}). Indeed, every blow-up limits must be $1$-noncollapsed by \cite{Brendle2}. If a blow-up limit is compact and strictly convex, then the original flow eventually becomes convex, and therefore converges to a family of shrinking round spheres (cf. \cite{Huisken1}). If a blow-up limit is noncompact and strictly convex, then that blow-up limit is the bowl soliton by Theorem \ref{main.thm}. Finally, if a blow-up limit is not strictly convex, then it splits off a line, and is a family of shrinking round cylinders.

We next discuss the background of Theorem \ref{main.thm}. Note that an ancient solution which is mean convex and two-sided noncollapsed is necessarily convex (cf. \cite{Haslhofer-Kleiner1}), but it is not uniformly convex unless it is a family of shrinking spheres (see \cite{Huisken-Sinestrari4}). Theorem \ref{main.thm} can be viewed as a parabolic analogue of the classical Bernstein theorem, which classifies entire solutions to the minimal surface equation. Theorem \ref{main.thm} can be generalized to higher dimensions, if we assume that the solution is uniformly two-convex (see \cite{Brendle-Choi}).

Daskalopoulos, Hamilton, and \v Se\v sum obtained a complete classification of all compact ancient solutions to the Ricci flow in dimension $2$ (cf. \cite{Daskalopoulos-Hamilton-Sesum2}). Moreover, they were able to classify all compact, convex ancient solutions to curve shortening flow in the plane (cf. \cite{Daskalopoulos-Hamilton-Sesum1}). Remarkably, these results do not require any noncollapsing assumptions. In a very important recent paper \cite{Angenent-Daskalopoulos-Sesum}, Angenent, Daskalopoulos, and \v Se\v sum studied compact, convex ancient solutions to the mean curvature flow. Under suitable symmetry assumptions, they obtained precise asymptotic estimates for the solution as $t \to -\infty$. 

A special case of ancient solutions are solitons; these are solutions that move in a self-similar fashion under the evolution. In a recent paper \cite{Brendle1}, the first author proved that every noncollapsed steady Ricci soliton in dimension $3$ is rotationally symmetric, and hence is isometric to the Bryant soliton up to scaling. Using similar techniques, Haslhofer \cite{Haslhofer} subsequently proved that every noncollapsed, convex translating soliton for the mean curvature flow in $\mathbb{R}^3$ is rotationally symmetric, and hence coincides with the bowl soliton up to scaling and ambient isometries. A related uniqueness result for the bowl soliton was proved in an important paper by Wang \cite{Wang}; this relies on a completely different approach. 

In Section \ref{asymptotic.analysis}, we study the asymptotic behavior of the flow as $t \to -\infty$. To that end, we write $M_t = (-t)^{\frac{1}{2}} \, \bar{M}_{-\log(-t)}$. As $\tau \to -\infty$, the rescaled surfaces $\bar{M}_\tau$ converge in $C_{\text{\rm loc}}^\infty$ to a cylinder of radius $\sqrt{2}$ with axis passing through the origin. More precisely, we show that $\bar{M}_\tau$ can be approximated by a cylinder up to error terms of order $O(e^{\frac{\tau}{2}})$. As in \cite{Colding-Minicozzi}, a major difficulty is the presence of a non-trivial eigenfunction for the linearized problem with eigenvalue $0$. This eigenfunction corresponds to the second Hermite polynomial. Using the convexity of $\bar{M}_\tau$ and the Brunn-Minkowski inequality, we show that this eigenfunction cannot become dominant as $\tau \to -\infty$.

In Section \ref{barrier}, we show that $\liminf_{t \to -\infty} H_{\text{\rm max}}(t) > 0$, where $H_{\text{\rm max}}(t)$ denotes the supremum of the mean curvature at time $t$. To do that, we consider the complement $M_t \setminus B_{8(-t)^{\frac{1}{2}}}(0)$. This set has two connected components, one of which is compact and one of which is noncompact. By combining the asymptotic analysis in Section \ref{asymptotic.analysis} with a barrier argument, we prove that the bounded connected component of $M_t \setminus B_{8(-t)^{\frac{1}{2}}}(0)$ has diameter at least $\sim (-t)$. This implies that $H_{\text{\rm max}}(t)$ is uniformly bounded from below as $t \to -\infty$. 

In Section \ref{neck.improvement}, we establish the Neck Improvement Theorem, which asserts that a neck becomes more symmetric under the evolution. This result does not require that the solution is ancient; it can be applied whenever we have a solution of mean curvature flow which is close to a cylinder on a sufficiently large parabolic neighborhood.

In Section \ref{rotational.symmetry}, we iterate the Neck Improvement Theorem to conclude that $M_t$ is rotationally symmetric. 

Finally, in Section \ref{analysis.in.rotationally.symmetric.case}, we analyze ancient solutions with rotational symmetry, and complete the proof of Theorem \ref{main.thm}.

\section{Asymptotic analysis as $t \to -\infty$}

\label{asymptotic.analysis}

Let $M_t$, $t \in (-\infty,0]$, be a noncompact ancient solution of mean curvature flow in $\mathbb{R}^3$ which is strictly convex and noncollapsed. We first recall some known results concerning the blowdown limit as $t \to -\infty$. Given any sequence $t_j \to -\infty$, we can find a subsequence with the property that the rescaled surfaces $(-t_j)^{-\frac{1}{2}} \, M_{t_j}$ converge in $C_{\text{\rm loc}}^\infty$ to a smooth limit, which is either a plane, or a round sphere, or a cylinder of radius $\sqrt{2}$ with axis passing through the origin (see \cite{Haslhofer-Kleiner1}, Theorem 1.11). Since the original flow $M_t$ is noncompact, the limit cannot be a sphere. Moreover, it follows from Huisken's monotonicity formula \cite{Huisken2} that the limit cannot be a plane. Therefore, the limit must be a cylinder. 

In the following, we consider the rescaled flow $\bar{M}_\tau = e^{\frac{\tau}{2}} \, M_{-e^{-\tau}}$. The surfaces $\bar{M}_\tau$ move with velocity $-(H - \frac{1}{2} \, \langle x,\nu \rangle) \nu$. Given any sequence $\tau_j \to -\infty$, there exists a subsequence with the property that the surfaces $\bar{M}_{\tau_j}$ converges in $C_{\text{\rm loc}}^\infty$ to a cylinder of radius $\sqrt{2}$ with axis passing through the origin. To fix notation, we denote by $\Sigma = \{x_1^2+x_2^2=2\}$ the cylinder of radius $\sqrt{2}$ whose axis coincides with the $x_3$-axis. 

\begin{proposition}
\label{gaussian.area}
For each $\tau$, we have 
\[\int_{\bar{M}_\tau} e^{-\frac{|x|^2}{4}} \leq \int_\Sigma e^{-\frac{|x|^2}{4}}.\]
\end{proposition}

\textbf{Proof.} 
Every convex surface is star-shaped, hence outward-minimizing by a standard calibration argument. This implies $\text{\rm area}(\bar{M}_\tau \cap B_r(p)) \leq Cr^2$ for all $p \in \mathbb{R}^3$ and all $r>0$. We next consider an arbitrary sequence $\tau_j \to -\infty$. After passing to a subsequence, the surfaces $\bar{M}_{\tau_j}$ converge in $C_{\text{\rm loc}}^\infty$ to a cylinder of radius $\sqrt{2}$ with axis passing through the origin. This gives 
\[\int_{\bar{M}_{\tau_j}} e^{-\frac{|x|^2}{4}} \to \int_\Sigma e^{-\frac{|x|^2}{4}}\] 
as $j \to \infty$. On the other hand, Huisken's monotonicity formula \cite{Huisken2} implies that the function 
\[\tau \mapsto \int_{\bar{M}_\tau} e^{-\frac{|x|^2}{4}}\] 
is monotone decreasing in $\tau$. From this, the assertion follows. \\

In view of the discussion above, there exists a smooth function $S(\tau)$ taking values in $SO(3)$ such that the rotated surfaces $\tilde{M}_\tau = S(\tau) \bar{M}_\tau$ converge to the cylinder $\Sigma$ in $C_{\text{\rm loc}}^\infty$. Hence, we can find a function $\rho(\tau)$ with the following properties: 
\begin{itemize} 
\item $\lim_{\tau \to -\infty} \rho(\tau) = \infty$. 
\item $-\rho(\tau) \leq \rho'(\tau) \leq 0$. 
\item We may write the surface $\tilde{M}_\tau$ as a graph of some function $u(\cdot,\tau)$ over $\Sigma \cap B_{2\rho(\tau)}(0)$, so that 
\[\{x + u(x,\tau) \nu_\Sigma(x): x \in \Sigma \cap B_{2\rho(\tau)}(0)\} \subset \tilde{M}_\tau,\] 
where $\nu_\Sigma$ denotes the unit normal to $\Sigma$ and $\|u(\cdot,\tau)\|_{C^4(\Sigma \cap B_{2\rho(\tau)}(0))} \leq \rho(\tau)^{-4}$.
\end{itemize} 
In the next step, we fine-tune the choice of $S(\tau)$. To that end, we fix a smooth cutoff function $\varphi \geq 0$ satisfying $\varphi=1$ on $[-\frac{1}{2},\frac{1}{2}]$ and $\varphi=0$ outside $[-\frac{2}{3},\frac{2}{3}]$. By the implicit function theorem, we can choose $S(\tau) \in SO(3)$ such that $u(x,\tau)$ satisfies the orthogonality relations 
\[\int_{\Sigma \cap B_{\rho(\tau)}(0)} e^{-\frac{|x|^2}{4}} \, \langle Ax,\nu_\Sigma \rangle \, u(x,\tau) \, \varphi \Big ( \frac{x_3}{\rho(\tau)} \Big ) = 0\] 
for every matrix $A \in so(3)$. Finally, we can arrange that the matrix $A(\tau) := S'(\tau) S(\tau)^{-1} \in so(3)$ satisfies $A(\tau)_{12}=0$. (Otherwise, we replace $S(\tau)$ by $R(\tau)S(\tau)$, where $R(\tau)$ is a rotation in the $x_1x_2$-plane. This does not affect the orthogonality relations above.) 

Our next two results are straightforward adaptations of powerful estimates in \cite{Angenent-Daskalopoulos-Sesum}. These estimates will play a key role in the subsequent arguments. Recall the foliation by self-shrinkers $\Sigma_a$ and $\tilde{\Sigma}_b$ given in \cite{Angenent-Daskalopoulos-Sesum} (see also \cite{Kleene-Moller}). Let $\nu_{\text{\rm fol}}$ denote the unit normal vector to this foliation. As explained in \cite{Angenent-Daskalopoulos-Sesum}, the union of all the leaves in this foliation contains a neighborhood of the cylinder $\Sigma$, and opens up like a cone at infinity. In particular, the union of all the leaves in the foliation contains the cylinder $\Sigma$ and the surface $\tilde{M}_\tau$ (for $-\tau$ sufficiently large), as well as the region $\Delta_\tau$ which lies in between $\Sigma$ and $\tilde{M}_\tau$. 

\begin{proposition}[cf. \cite{Angenent-Daskalopoulos-Sesum}, Lemma 4.10]
\label{calibration}
There exists a constant $L_0$ such that for all $L \in [L_0,\rho(\tau)]$ 
\[\int_{\tilde{M}_\tau \cap \{|x_3| \geq L\}} e^{-\frac{|x|^2}{4}} - \int_{\Sigma \cap \{|x_3| \geq L\}} e^{-\frac{|x|^2}{4}} \geq -\int_{\Delta_\tau \cap \{|x_3|=L\}} e^{-\frac{x^2}{4}} \, |\langle \omega,\nu_{\text{\rm fol}} \rangle|,\] 
where $\omega=(0,0,1)$ denotes the vertical unit vector in $\mathbb{R}^3$.
\end{proposition}

\textbf{Proof.} 
Since each leaf of the foliation is a self-similar shrinker, the normal vector $\nu_{\text{\rm fol}}$ satisfies $\text{\rm div}(e^{-\frac{|x|^2}{4}} \, \nu_{\text{\rm fol}}) = 0$. Note that $\nu_{\text{\rm fol}} = \nu_\Sigma$ at each point on the cylinder $\Sigma$. Using the divergence theorem, we obtain 
\begin{align*} 
&\int_{\tilde{M}_\tau \cap \{L \leq |x_3| \leq z\}} e^{-\frac{|x|^2}{4}} \, \langle \nu,\nu_{\text{\rm fol}} \rangle - \int_{\Sigma \cap \{L \leq |x_3| \leq z\}} e^{-\frac{|x|^2}{4}} \\ 
&\geq -\int_{\Delta_\tau \cap \{|x_3|=L\}} e^{-\frac{x^2}{4}} \, |\langle \omega,\nu_{\text{\rm fol}} \rangle| - \int_{\Delta_\tau \cap \{|x_3|=z\}} e^{-\frac{x^2}{4}} \, |\langle \omega,\nu_{\text{\rm fol}} \rangle|. 
\end{align*} 
We know $\langle \nu,\nu_{\text{\rm fol}} \rangle \leq 1$ on $\tilde{M}_\tau$. The convexity of $\tilde{M}_\tau$ implies that the area of $\Delta_\tau \cap \{|x_3|=z\}$ is at most $Cz^2$. Hence, passing to the limit as $z \to \infty$ gives 
\[\int_{\tilde{M}_\tau \cap \{|x_3| \geq L\}} e^{-\frac{|x|^2}{4}} - \int_{\Sigma \cap \{|x_3| \geq L\}} e^{-\frac{|x|^2}{4}} \geq -\int_{\Delta_\tau \cap \{|x_3|=L\}} e^{-\frac{x^2}{4}} \, |\langle \omega,\nu_{\text{\rm fol}} \rangle|,\] 
as desired. \\

\begin{proposition}[cf. \cite{Angenent-Daskalopoulos-Sesum}, Lemma 4.7]
\label{consequence.of.gaussian.area.bound}
There exists a constant $L_0$ such that 
\[\int_{\Sigma \cap \{|x_3| \leq L\}} e^{-\frac{|x|^2}{4}} \, |\nabla^\Sigma u(x,\tau)|^2 \leq C \int_{\Sigma \cap \{|x_3| \leq \frac{L}{2}\}} e^{-\frac{|x|^2}{4}} \, u(x,\tau)^2\] 
and 
\[\int_{\Sigma \cap \{\frac{L}{2} \leq |x_3| \leq L\}} e^{-\frac{|x|^2}{4}} \, u(x,\tau)^2 \leq CL^{-2} \int_{\Sigma \cap \{|x_3| \leq \frac{L}{2}\}} e^{-\frac{|x|^2}{4}} \, u(x,\tau)^2\] 
for all $L \in [L_0,\rho(\tau)]$.
\end{proposition}

\textbf{Proof.} 
Lemma 4.11 in \cite{Angenent-Daskalopoulos-Sesum} implies that $|\langle \omega,\nu_{\text{\rm fol}} \rangle| \leq C L^{-1} \, |x_1^2+x_2^2-2|$ for each point $x \in \Delta_\tau \cap \{|x_3|=L\}$. This gives 
\begin{align*} 
\int_{\Delta_\tau \cap \{|x_3|=L\}} e^{-\frac{x^2}{4}} \, |\langle \omega,\nu_{\text{\rm fol}} \rangle| 
&\leq CL^{-1} \int_{\Delta_\tau \cap \{|x_3|=L\}} e^{-\frac{x^2}{4}} \, |x_1^2+x_2^2-2| \\ 
&\leq C L^{-1} \int_{\Sigma \cap \{|x_3|=L\}} e^{-\frac{|x|^2}{4}} \, u^2. 
\end{align*}
Using Proposition \ref{calibration}, we obtain 
\[\int_{\tilde{M}_\tau \cap \{|x_3| \geq L\}} e^{-\frac{|x|^2}{4}} - \int_{\Sigma \cap \{|x_3| \geq L\}} e^{-\frac{|x|^2}{4}} \geq -C L^{-1} \int_{\Sigma \cap \{|x_3|=L\}} e^{-\frac{|x|^2}{4}} \, u^2\] 
(compare \cite{Angenent-Daskalopoulos-Sesum}, equation (4.33)). On the other hand, 
\begin{align*} 
&\int_{\tilde{M}_\tau \cap \{|x_3| \leq L\}} e^{-\frac{|x|^2}{4}} - \int_{\Sigma \cap \{|x_3| \leq L\}} e^{-\frac{|x|^2}{4}} \\ 
&= \int_{-L}^L \int_0^{2\pi} e^{-\frac{z^2}{4}} \, \bigg [ e^{-\frac{(\sqrt{2}+u)^2}{4}} \, \sqrt{(\sqrt{2}+u)^2 \, \Big ( 1+ \Big ( \frac{\partial u}{\partial z} \Big )^2 \Big ) + \Big ( \frac{\partial u}{\partial \theta} \Big )^2} - e^{-\frac{1}{2}} \, \sqrt{2} \bigg ] \, d\theta \, dz.
\end{align*}
Since $L \leq \rho(\tau)$, we have $|u|+|\frac{\partial u}{\partial z}|+|\frac{\partial u}{\partial \theta}| \leq o(1)$ for $|x_3| \leq L$. This gives 
\begin{align*} 
&\int_{\tilde{M}_\tau \cap \{|x_3| \leq L\}} e^{-\frac{|x|^2}{4}} - \int_{\Sigma \cap \{|x_3| \leq L\}} e^{-\frac{|x|^2}{4}} \\ 
&\geq \int_{-L}^L \int_0^{2\pi} e^{-\frac{z^2}{4}} \, \Big [ e^{-\frac{(\sqrt{2}+u)^2}{4}} \, (\sqrt{2}+u) - e^{-\frac{1}{2}} \, \sqrt{2} + \frac{1}{C} \, |\nabla^\Sigma u|^2 \Big ] \, d\theta \, dz \\ 
&\geq \int_{-L}^L \int_0^{2\pi} e^{-\frac{z^2}{4}} \, \Big [ -C u^2 + \frac{1}{C} \, |\nabla^\Sigma u|^2 \Big ] \, d\theta \, dz
\end{align*} 
where $C>0$ is a large numerical constant. Putting these facts together, we obtain 
\begin{align*} 
\int_{\tilde{M}_\tau} e^{-\frac{|x|^2}{4}} - \int_\Sigma e^{-\frac{|x|^2}{4}} 
&\geq \int_{\Sigma \cap \{|x_3| \leq L\}} e^{-\frac{z^2}{4}} \, \Big [ -C u^2 + \frac{1}{C} \, |\nabla^\Sigma u|^2 \Big ] \\ 
&- C L^{-1} \int_{\Sigma \cap \{|x_3|=L\}} e^{-\frac{|x|^2}{4}} \, u^2.
\end{align*} 
Using Proposition \ref{gaussian.area}, we conclude that 
\begin{align*} 
&\int_{\Sigma \cap \{|x_3| \leq L\}} e^{-\frac{|x|^2}{4}} \, |\nabla^\Sigma u|^2 \\
&\leq C \int_{\Sigma \cap \{|x_3| \leq L\}} e^{-\frac{|x|^2}{4}} u^2 + C L^{-1} \int_{\Sigma \cap \{|x_3|=L\}} e^{-\frac{|x|^2}{4}} \, u^2. 
\end{align*}
The divergence theorem gives 
\begin{align*} 
L \int_{\Sigma \cap \{|x_3|=L\}} e^{-\frac{|x|^2}{4}} \, u^2 
&= \int_{\Sigma \cap \{|x_3| \leq L\}} \text{\rm div}_\Sigma(e^{-\frac{|x|^2}{4}} \, u^2 \, x^{\text{\rm tan}}) \\ 
&= \int_{\Sigma \cap \{|x_3| \leq L\}} e^{-\frac{|x|^2}{4}} \, \Big ( u^2 - \frac{1}{2} \, x_3^2 \, u^2 + 2u \, \langle x^{\text{\rm tan}},\nabla^\Sigma u \rangle \Big ) \\ 
&\leq \int_{\Sigma \cap \{|x_3| \leq L\}} e^{-\frac{|x|^2}{4}} \, \Big ( u^2 - \frac{1}{4} \, x_3^2 \, u^2 + 4 \, |\nabla^\Sigma u|^2 \Big ),
\end{align*}
hence 
\begin{align*} 
&L^2 \int_{\Sigma \cap \{|x_3| \leq L\}} e^{-\frac{|x|^2}{4}} u^2 + L \int_{\Sigma \cap \{|x_3|=L\}} e^{-\frac{|x|^2}{4}} \, u^2 \\ 
&\leq C \int_{\Sigma \cap \{|x_3| \leq L\}} e^{-\frac{|x|^2}{4}} \, |\nabla^\Sigma u|^2 + C L^2 \int_{\Sigma \cap \{|x_3| \leq \frac{L}{2}\}} e^{-\frac{|x|^2}{4}} \, u^2. 
\end{align*}
Putting these facts together, we conclude that 
\begin{align*} 
&\int_{\Sigma \cap \{|x_3| \leq L\}} e^{-\frac{|x|^2}{4}} \, |\nabla^\Sigma u|^2 \\
&\leq CL^{-2} \int_{\Sigma \cap \{|x_3| \leq L\}} e^{-\frac{|x|^2}{4}} \, |\nabla^\Sigma u|^2 + C \int_{\Sigma \cap \{|x_3| \leq \frac{L}{2}\}} e^{-\frac{|x|^2}{4}} \, u^2. 
\end{align*} 
If $L$ is sufficiently large, the first term on the right hand side can be absorbed into the left hand side. This gives  
\[\int_{\Sigma \cap \{|x_3| \leq L\}} e^{-\frac{|x|^2}{4}} \, |\nabla^\Sigma u|^2 \leq C \int_{\Sigma \cap \{|x_3| \leq \frac{L}{2}\}} e^{-\frac{|x|^2}{4}} \, u^2.\] 
This proves the first statement. Using the inequality 
\[0 \leq \int_{\Sigma \cap \{|x_3| \leq L\}} e^{-\frac{|x|^2}{4}} \, \Big ( u^2 - \frac{1}{4} \, x_3^2 \, u^2 + 4 \, |\nabla^\Sigma u|^2 \Big ),\] 
the second statement follows. \\

Let us denote by $\mathcal{H}$ the space of all functions $f$ on $\Sigma$ such that 
\[\|f\|_{\mathcal{H}}^2 = \int_\Sigma e^{-\frac{|x|^2}{4}} \, f^2 < \infty.\] 
We define an operator $\mathcal{L}$ on the cylinder $\Sigma$ by 
\[\mathcal{L} f = \Delta_\Sigma f - \frac{1}{2} \, \langle x^{\text{\rm tan}},\nabla^\Sigma f \rangle + f.\] 
In coordinates, $\mathcal{L}$ takes the form 
\[\mathcal{L} f = \frac{\partial^2}{\partial z^2} f + \frac{1}{2} \, \frac{\partial^2}{\partial \theta^2} f - \frac{1}{2} \, z \, \frac{\partial}{\partial z} f + f.\] 
The eigenfunctions of $\mathcal{L}$ are of the form $H_n \big ( \frac{z}{2} \big ) \, \cos(m\theta)$ and $H_n \big ( \frac{z}{2} \big ) \, \sin(m\theta)$, where $m$ and $n$ are nonnegative integers and $H_n$ denotes the Hermite polynomial of degree $n$. The corresponding eigenvalues are given by $1-\frac{n+m^2}{2}$. Thus, there are four eigenfunctions that correspond to positive eigenvalues of $\mathcal{L}$, and these are given by $1$, $z$, $\cos \theta$, $\sin \theta$ up to scaling. The span of these eigenfunctions will be denoted by $\mathcal{H}_+$. Moreover, there are three eigenfunctions of $\mathcal{L}$ with eigenvalue $0$, and these are given by $z^2-2$, $z \cos \theta$, $z \sin \theta$ up to scaling. The span of these eigenfunctions will be denoted by $\mathcal{H}_0$. The span of all other eigenfunctions will be denoted by $\mathcal{H}_-$. Clearly, 
\begin{align*} 
&\langle \mathcal{L} f,f \rangle_{\mathcal{H}} \geq \frac{1}{2} \, \|f\|_{\mathcal{H}}^2 & \text{\rm for $f \in \mathcal{H}_+$,} \\ 
&\langle \mathcal{L} f,f \rangle_{\mathcal{H}} = 0 & \text{\rm for $f \in \mathcal{H}_0$,} \\ 
&\langle \mathcal{L} f,f \rangle_{\mathcal{H}} \leq -\frac{1}{2} \, \|f\|_{\mathcal{H}}^2 & \text{\rm for $f \in \mathcal{H}_-$.} 
\end{align*}

\begin{lemma}
\label{pde.for.u}
The function $u(x,\tau)$ satisfies 
\[\frac{\partial}{\partial \tau} u = \mathcal{L} u + E + \langle A(\tau)x,\nu_\Sigma \rangle,\] 
where $E$ satisfies the pointwise estimate $|E| \leq O(\rho(\tau)^{-1}) \, (|u| + |\nabla^\Sigma u| + |A(\tau)|)$.
\end{lemma}

\textbf{Proof.} 
Recall that the rescaled surfaces $\bar{M}_\tau$ move with velocity $-(H - \frac{1}{2} \, \langle x,\nu \rangle) \nu$. Hence, the rotated surfaces $\tilde{M}_\tau = S(\tau) \bar{M}_\tau$ move with velocity $-(H - \frac{1}{2} \, \langle x,\nu \rangle - \langle A(\tau)x,\nu \rangle) \nu$, where $x \in \tilde{M}_\tau$. Therefore, the function $u(x,\tau)$ satisfies the equation 
\begin{align*} 
\frac{\partial}{\partial \tau} u 
&= -\frac{1}{\langle \nu_\Sigma,\nu(x+u\nu_\Sigma) \rangle} \, H(x+u\nu_\Sigma) \\ 
&+ \frac{1}{2 \langle \nu_\Sigma,\nu(x+u\nu_\Sigma) \rangle} \, \langle x+u\nu_\Sigma,\nu(x+u\nu_\Sigma) \rangle \\ 
&+ \frac{1}{\langle \nu_\Sigma,\nu(x+u\nu_\Sigma) \rangle} \, \langle A(\tau)(x+u\nu_\Sigma),\nu(x+u\nu_\Sigma) \rangle) 
\end{align*}
for $x \in \Sigma$. By assumption, $\|u\|_{C^4(\Sigma \cap B_{2\rho(\tau)}(0))} \leq O(\rho(\tau)^{-2})$. This gives 
\[\big | \nu(x+u\nu_\Sigma) - \nu_\Sigma + \nabla^\Sigma u \big | \leq O(\rho(\tau)^{-2}) \, (|u| + |\nabla^\Sigma u|)\] 
and 
\[\Big | H(x+u\nu_\Sigma) + \Delta_\Sigma u + \frac{1}{2} \, u \Big | \leq O(\rho(\tau)^{-2}) \, (|u| + |\nabla^\Sigma u|).\] 
Putting these facts together, we obtain 
\[\frac{\partial}{\partial \tau} u = \mathcal{L} u + E + \langle A(\tau)x,\nu_\Sigma \rangle,\] 
where $E$ satisfies the pointwise estimate $|E| \leq O(\rho(\tau)^{-1}) \, (|u| + |\nabla^\Sigma u| + |A(\tau)|)$. \\

\begin{lemma} 
\label{pde.for.hat.u}
The function $\hat{u}(x,\tau) = u(x,\tau) \, \varphi \big ( \frac{x_3}{\rho(\tau)} \big )$ satisfies 
\[\frac{\partial}{\partial \tau} \hat{u} = \mathcal{L} \hat{u} + \hat{E} + \langle A(\tau)x,\nu_\Sigma \rangle \, \varphi \Big ( \frac{x_3}{\rho(\tau)} \Big ),\] 
where $\hat{E}$ satisfies $\|\hat{E}\|_{\mathcal{H}} \leq O(\rho(\tau)^{-1}) \, (\|\hat{u}\|_{\mathcal{H}} + |A(\tau)|)$.
\end{lemma} 

\textbf{Proof.} 
We compute 
\[\frac{\partial}{\partial \tau} \hat{u} = \mathcal{L} \hat{u} + \hat{E} + \langle A(\tau)x,\nu_\Sigma \rangle \, \varphi \Big ( \frac{x_3}{\rho(\tau)} \Big )\] 
where 
\begin{align*} 
\hat{E} &= E \, \varphi \Big ( \frac{x_3}{\rho(\tau)} \Big ) - \frac{2}{\rho(\tau)} \, \frac{\partial u}{\partial z} \, \varphi' \Big ( \frac{x_3}{\rho(\tau)} \Big ) - \frac{1}{\rho(\tau)^2} \, u \, \varphi'' \Big ( \frac{x_3}{\rho(\tau)} \Big ) \\ 
&+ \frac{x_3}{2\rho(\tau)} \, u \, \varphi' \Big ( \frac{x_3}{\rho(\tau)} \Big ) - \frac{x_3 \rho'(\tau)}{\rho(\tau)^2} \, u \, \varphi' \Big ( \frac{x_3}{\rho(\tau)} \Big ). 
\end{align*} 
Using Lemma \ref{pde.for.u}, we deduce that 
\[|\hat{E}| \leq O(\rho(\tau)^{-1}) \, (|u| + |\nabla^\Sigma u| + |A(\tau)|)\] 
for $|x_3| \leq \frac{\rho(\tau)}{2}$. Moreover, since $|\rho'(\tau)| \leq \rho(\tau)$, we obtain 
\[|\hat{E}| \leq O(1) \, |u| + O(\rho(\tau)^{-1}) \, (|\nabla^\Sigma u| + |A(\tau)|)\] 
for $\frac{\rho(\tau)}{2} \leq |x_3| \leq \rho(\tau)$. Using Proposition \ref{consequence.of.gaussian.area.bound}, we conclude that 
\begin{align*} 
\int_\Sigma e^{-\frac{|x|^2}{4}} \, |\hat{E}|^2 
&\leq O(\rho(\tau)^{-2}) \int_{\Sigma \cap \{|x_3| \leq \frac{\rho(\tau)}{2}\}} e^{-\frac{|x|^2}{4}} \, u^2 \\ 
&+ O(1) \int_{\Sigma \cap \{\frac{\rho(\tau)}{2} \leq |x_3| \leq \rho(\tau)\}} e^{-\frac{|x|^2}{4}} \, u^2 \\ 
&+ O(\rho(\tau)^{-2}) \int_{\Sigma \cap \{|x_3| \leq \rho(\tau)\}} e^{-\frac{|x|^2}{4}} \, |\nabla^\Sigma u|^2 \\ 
&+ O(\rho(\tau)^{-2}) \, |A(\tau)|^2 \\ 
&\leq O(\rho(\tau)^{-2}) \int_{\Sigma \cap \{|x_3| \leq \frac{\rho(\tau)}{2}\}} e^{-\frac{|x|^2}{4}} \, u^2 \\ 
&+ O(\rho(\tau)^{-2}) \, |A(\tau)|^2 \\ 
&\leq O(\rho(\tau)^{-2}) \int_\Sigma e^{-\frac{|x|^2}{4}} \, \hat{u}^2 \\ 
&+ O(\rho(\tau)^{-2}) \, |A(\tau)|^2. 
\end{align*}
Thus, $\|\hat{E}\|_{\mathcal{H}} \leq O(\rho(\tau)^{-1}) \, \|\hat{u}\|_{\mathcal{H}} + O(\rho(\tau)^{-1}) \, |A(\tau)|$, as claimed. \\

\begin{lemma} 
\label{pde.for.hat.u.2}
We have $|A(\tau)| \leq O(\rho(\tau)^{-1}) \, \|u\|_{\mathcal{H}}$ and 
\[\Big \| \frac{\partial}{\partial \tau} \hat{u} - \mathcal{L} \hat{u} \Big \|_{\mathcal{H}} \leq O(\rho(\tau)^{-1}) \, \|\hat{u}\|_{\mathcal{H}}.\] 
\end{lemma} 

\textbf{Proof.} 
The orthogonality relations imply that $\hat{u}$ is orthogonal to $\langle Ax,\nu_\Sigma \rangle$ for every $A \in so(3)$. Since this is true for each $\tau$, it follows that $\frac{\partial}{\partial \tau} \hat{u}$ is orthogonal to $\langle Ax,\nu_\Sigma \rangle$ for every $A \in so(3)$. Moreover, since the function $\langle Ax,\nu_\Sigma \rangle$ is an eigenfunction of $\mathcal{L}$ with eigenvalue $0$, we deduce that $\mathcal{L} \hat{u}$ is orthogonal to $\langle Ax,\nu_\Sigma \rangle$ for every $A \in so(3)$. Consequently, $\frac{\partial}{\partial \tau} \hat{u} - \mathcal{L} \hat{u}$ is orthogonal to $\langle Ax,\nu_\Sigma \rangle$ for every $A \in so(3)$. Therefore, $\hat{E} + \langle A(\tau)x,\nu_\Sigma \rangle \, \varphi \big ( \frac{x_3}{\rho(\tau)} \big )$ is orthogonal to $\langle Ax,\nu_\Sigma \rangle$ for every $A \in so(3)$. In particular, 
\[\int_\Sigma  e^{-\frac{|x|^2}{4}} \, \Big ( \hat{E} + \langle A(\tau)x,\nu_\Sigma \rangle \, \varphi \big ( \frac{x_3}{\rho(\tau)} \big ) \Big ) \, \langle A(\tau)x,\nu_\Sigma \rangle = 0.\] 
Using the fact that $A(\tau)_{12}=0$, we obtain 
\begin{align*} 
|A(\tau)|^2 
&\leq O(1) \int_\Sigma e^{-\frac{|x|^2}{4}} \, \langle A(\tau)x,\nu_\Sigma \rangle^2 \, \varphi \big ( \frac{x_3}{\rho(\tau)} \big ) \\ 
&\leq O(1) \int_\Sigma e^{-\frac{|x|^2}{4}} \, |\hat{E}| \, |\langle A(\tau)x,\nu_\Sigma \rangle| \\ 
&\leq O(1) \, \|\hat{E}\|_{\mathcal{H}} \, |A(\tau)| \\ 
&\leq O(\rho(\tau)^{-1}) \, (\|\hat{u}\|_{\mathcal{H}}+|A(\tau)|) \, |A(\tau)|, 
\end{align*} 
where in the last step we have used Lemma \ref{pde.for.hat.u}. Consequently, $|A(\tau)| \leq O(\rho(\tau)^{-1}) \, \|\hat{u}\|_{\mathcal{H}}$. Using Lemma \ref{pde.for.hat.u}, we obtain \[\Big \| \frac{\partial}{\partial \tau} \hat{u} - \mathcal{L} \hat{u} \Big \|_{\mathcal{H}} \leq \|\hat{E}\|_{\mathcal{H}} + O(1) \, |A(\tau)| \leq O(\rho(\tau)^{-1}) \, \|\hat{u}\|_{\mathcal{H}},\] 
as claimed. \\

We now define 
\begin{align*} 
&U_+(\tau) := \|P_+ \hat{u}(\cdot,\tau)\|_{\mathcal{H}}^2, \\ 
&U_0(\tau) := \|P_0 \hat{u}(\cdot,\tau)\|_{\mathcal{H}}^2, \\ 
&U_-(\tau) := \|P_- \hat{u}(\cdot,\tau)\|_{\mathcal{H}}^2, 
\end{align*} 
where $P_+, P_0, P_-$ denote the orthogonal projections to $\mathcal{H}_+,\mathcal{H}_0,\mathcal{H}_-$, respectively. Using Lemma \ref{pde.for.hat.u.2}, we obtain 
\begin{align*} 
&\frac{d}{d\tau} U_+(\tau) \geq U_+(\tau) - O(\rho(\tau)^{-1}) \, (U_+(\tau) + U_0(\tau) + U_-(\tau)), \\ 
&\Big | \frac{d}{d\tau} U_0(\tau) \Big | \leq O(\rho(\tau)^{-1}) \, (U_+(\tau) + U_0(\tau) + U_-(\tau)), \\ 
&\frac{d}{d\tau} U_-(\tau) \leq -U_-(\tau) + O(\rho(\tau)^{-1}) \, (U_+(\tau) + U_0(\tau) + U_-(\tau)). 
\end{align*}
Clearly, $U_+(\tau)+U_0(\tau)+U_-(\tau) = \|\hat{u}\|_{\mathcal{H}}^2 \to 0$ as $\tau \to -\infty$. Moreover, $U_+(\tau)+U_0(\tau)+U_-(\tau) = \|\hat{u}\|_{\mathcal{H}}^2 > 0$ since $\tilde{M}_\tau$ is strictly convex. 

\begin{lemma}
\label{consequence.of.merle.zaag}
We have $U_0(\tau)+U_-(\tau) \leq o(1) U_+(\tau)$.
\end{lemma} 

\textbf{Proof.} 
Applying an ODE lemma of Merle and Zaag (cf. Lemma 5.4 in \cite{Angenent-Daskalopoulos-Sesum} or Lemma A.1 in \cite{Merle-Zaag}), we conclude that either $U_0(\tau)+U_-(\tau) \leq o(1) U_+(\tau)$ or $U_+(\tau)+U_-(\tau) \leq o(1) U_0(\tau)$.

The second case can be ruled out as follows: Suppose $U_+(\tau)+U_-(\tau) \leq o(1) U_0(\tau)$. Then $\frac{\hat{u}(\cdot,\tau)}{\|\hat{u}(\cdot,\tau)\|_{\mathcal{H}}}$ converges with respect to $\|\cdot\|_{\mathcal{H}}$ to the subspace $\mathcal{H}_0 = \text{\rm span}\{z^2-2,z \cos \theta,z \sin \theta\}$. On the other hand, the orthogonality relations used to define $S(\tau)$ imply that the function $\hat{u}(\cdot,\tau)$ is orthogonal to the function $\langle Ax,\nu_\Sigma \rangle$ for each $A \in so(3)$. In other words, the function $\hat{u}(\cdot,\tau)$ is orthogonal to the functions $z \cos \theta$ and $z \sin \theta$. Consequently, $\frac{\hat{u}(\cdot,\tau)}{\|\hat{u}(\cdot,\tau)\|_{\mathcal{H}}}$ converges with respect to $\|\cdot\|_{\mathcal{H}}$ to a non-zero multiple of $z^2-2$. 

Let $\Omega_\tau$ denote the region enclosed by $\tilde{M}_\tau$. Moreover, we denote by $\mathcal{A}(z,\tau)$ the area of the intersection $\Omega_\tau \cap \{x_3=z\}$. By the Brunn-Minkowski inequality, the function $z \mapsto \sqrt{\mathcal{A}(z,\tau)}$ is concave. Since $\tilde{M}_\tau$ is noncompact, it follows that the function $z \mapsto \sqrt{\mathcal{A}(z,\tau)}$ is monotone. 

For $|z| \leq \rho(\tau)$, we have the exact identity $\mathcal{A}(z,\tau) = \frac{1}{2} \int_0^{2\pi} (\sqrt{2}+u(\theta,z,\tau))^2 \, d\theta$. Thus, the function $z \mapsto \int_0^{2\pi} (2\sqrt{2} \, u(\theta,z,\tau)+u(\theta,z,\tau)^2) \, d\theta$ is monotone. In particular, we either have 
\begin{align*} 
&\int_{-3}^{-1} \int_0^{2\pi} (2\sqrt{2} \, u(\theta,z,\tau)+u(\theta,z,\tau)^2) \, d\theta \, dz \\ 
&\leq \int_{-1}^1 \int_0^{2\pi} (2\sqrt{2} \, u(\theta,z,\tau)+u(\theta,z,\tau)^2) \, d\theta \, dz \\ 
&\leq \int_1^3 \int_0^{2\pi} (2\sqrt{2} \, u(\theta,z,\tau)+u(\theta,z,\tau)^2) \, d\theta \, dz 
\end{align*} 
or 
\begin{align*} 
&\int_{-3}^{-1} \int_0^{2\pi} (2\sqrt{2} \, u(\theta,z,\tau)+u(\theta,z,\tau)^2) \, d\theta \, dz \\ 
&\geq \int_{-1}^1 \int_0^{2\pi} (2\sqrt{2} \, u(\theta,z,\tau)+u(\theta,z,\tau)^2) \, d\theta \, dz \\ 
&\geq \int_1^3 \int_0^{2\pi} (2\sqrt{2} \, u(\theta,z,\tau)+u(\theta,z,\tau)^2) \, d\theta \, dz. 
\end{align*} 
However, neither of these possibilities is consistent with the fact that the limit of $\frac{\hat{u}(\cdot,\tau)}{\|\hat{u}(\cdot,\tau)\|_{\mathcal{H}}}$ is a non-zero multiple of $z^2-2$. This is a contradiction. This completes the proof of Lemma \ref{consequence.of.merle.zaag}. \\

\begin{lemma} 
\label{almost.sharp.asymptotics.for.u}
For each $\varepsilon > 0$, we have $\|u(\cdot,\tau)\|_{C^4([0,2\pi] \times [-100,100])} \leq o(e^{\frac{(1-\varepsilon)\tau}{2}})$ and $|A(\tau)| \leq o(e^{\frac{(1-\varepsilon)\tau}{2}})$.
\end{lemma}

\textbf{Proof.} 
Recall that $U_0(\tau)+U_-(\tau) \leq o(1) U_+(\tau)$ by Lemma \ref{consequence.of.merle.zaag}. This directly implies 
\[\frac{d}{d\tau} U_+(\tau) \geq U_+(\tau) - o(1) \, U_+(\tau).\] 
Consequently, for every $\varepsilon>0$, we have $U_+(\tau) \leq o(e^{(1-\varepsilon)\tau})$. This gives $U_0(\tau)+U_-(\tau) \leq o(1) U_+(\tau) \leq o(e^{(1-\varepsilon)\tau})$, hence 
\[\|\hat{u}\|_{\mathcal{H}}^2 = U_+(\tau) + U_0(\tau) + U_-(\tau) \leq o(e^{(1-\varepsilon)\tau}).\] 
Using Lemma \ref{pde.for.hat.u.2}, we obtain $|A(\tau)| \leq o(1) \, \|\hat{u}\|_{\mathcal{H}} \leq o(e^{\frac{(1-\varepsilon)\tau}{2}})$. Finally, the inequality $\|u(\cdot,\tau)\|_{C^4([0,2\pi] \times [-100,100])} \leq o(e^{\frac{(1-\varepsilon)\tau}{2}})$ follows from standard interpolation inequalities. This completes the proof of Lemma \ref{almost.sharp.asymptotics.for.u}. \\

Recall that $A(\tau) = S'(\tau) S(\tau)^{-1}$. Since $|A(\tau)| \leq o(e^{\frac{(1-\varepsilon)\tau}{2}})$ by Lemma \ref{almost.sharp.asymptotics.for.u}, the limit $\lim_{\tau \to -\infty} S(\tau)$ exists. Without loss of generality, we may assume that $\lim_{\tau \to -\infty} S(\tau) = \text{\rm id}$. Clearly, $|S(\tau)-\text{\rm id}| \leq o(e^{\frac{(1-\varepsilon)\tau}{2}})$.

\begin{lemma} 
\label{asymptotics.for.surface}
We have 
\[\sup_{\bar{M}_\tau \cap \{|x_3| \leq e^{-\frac{\tau}{10}}\}} |x_1^2+x_2^2-2| \leq e^{\frac{\tau}{10}}\] 
if $-\tau$ is sufficiently large.
\end{lemma}

\textbf{Proof.} 
Using Lemma \ref{almost.sharp.asymptotics.for.u} and the estimate $|S(\tau)-\text{\rm id}| \leq o(e^{\frac{(1-\varepsilon)\tau}{2}})$, we obtain 
\[\sup_{x \in \bar{M}_\tau \cap B_{10}(0)} |x_1^2+x_2^2-2| \leq o(e^{\frac{(1-\varepsilon)\tau}{2}}).\] 
The convexity of $\bar{M}_\tau$ implies 
\[\sup_{\bar{M}_\tau \cap \{|x_3| \leq e^{-\frac{\tau}{10}}\}} (x_1^2+x_2^2) \leq 2 + e^{\frac{\tau}{10}}\] 
if $-\tau$ is sufficiently large. Let 
\[\Sigma_a = \{(x_1,x_2,x_3): x_1^2+x_2^2=u_a(-x_3)^2, \, -a \leq x_3 \leq 0\}\] 
denote the self-similar shrinker constructed in \cite{Angenent-Daskalopoulos-Sesum}. By Lemma 4.4 in \cite{Angenent-Daskalopoulos-Sesum}, $u_a(2) \leq \sqrt{2}-a^{-2}$. Since $\bar{M}_\tau$ converges to $\Sigma$ in $C_{\text{\rm loc}}^\infty$, the surface $\bar{M}_\tau \cap \{x_3 \leq -2\}$ encloses the surface $\Sigma_a \cap \{x_3 \leq -2\}$ if $-\tau$ is sufficiently large (depending on $a$). On the other hand, since $\inf_{x \in \bar{M}_\tau \cap B_{10}(0)} (x_1^2+x_2^2) \geq 2 - o(e^{\frac{(1-\varepsilon)\tau}{2}})$, the boundary $\bar{M}_\tau \cap \{x_3=-2\}$ encloses the boundary $\Sigma_a \cap \{x_3=-2\}$ provided that $-\tau$ is sufficiently large and $a \leq e^{-\frac{(1-\varepsilon)\tau}{4}}$. By the maximum principle, the surface $\bar{M}_\tau \cap \{x_3 \leq -2\}$ encloses $\Sigma_a \cap \{x_3 \leq -2\}$ whenever $-\tau$ is sufficiently large and $a \leq e^{-\frac{(1-\varepsilon)\tau}{4}}$. Using Theorem 8.2 in \cite{Angenent-Daskalopoulos-Sesum}, we obtain $u_a(y) \geq \sqrt{2(1-a^{-2}y^2)}$ for all $y \in [0,a]$, provided that $a$ is sufficiently large. Putting these facts together, we obtain 
\[\inf_{\bar{M}_\tau \cap \{-e^{-\frac{\tau}{10}} \leq x_3 \leq -2\}} (x_1^2+x_2^2) \geq 2 - e^{\frac{\tau}{10}}\] 
if $-\tau$ is sufficiently large. An analogous argument gives 
\[\inf_{\bar{M}_\tau \cap \{2 \leq x_3 \leq e^{-\frac{\tau}{10}}\}} (x_1^2+x_2^2) \geq 2 - e^{\frac{\tau}{10}}\] 
if $-\tau$ is sufficiently large. Putting these facts together, we obtain 
\[\inf_{\bar{M}_\tau \cap \{|x_3| \leq e^{-\frac{\tau}{10}}\}} (x_1^2+x_2^2) \geq 2-e^{\frac{\tau}{10}}\] 
if $-\tau$ is sufficiently large. This completes the proof of Lemma \ref{asymptotics.for.surface}. \\

\begin{lemma}
\label{size.of.necklike.region} 
Let $\varepsilon_0>0$ be given. If $-\tau$ is sufficiently large (depending on $\varepsilon_0$), then every point in $\bar{M}_\tau \cap \{|x_3| \leq \frac{1}{2} \, e^{-\frac{\tau}{10}}\}$ lies at the center of an $\varepsilon_0$-neck.
\end{lemma}

\textbf{Proof.} 
Suppose that there exists a sequence of times $\tau_j \to -\infty$ and a sequence of points $q_j \in \bar{M}_{\tau_j} \cap \{|x_3| \leq \frac{1}{2} \, e^{-\frac{\tau_j}{10}}\}$ such that $q_j$ does not lie on an $\varepsilon_0$-neck. Using Lemma \ref{asymptotics.for.surface} and the noncollapsing property, we conclude that the mean curvature at $q_j$ is bounded from below by a positive constant. We now consider the triangle in $\mathbb{R}^3$ with vertices $q_j$, $(0,0,e^{-\frac{\tau_j}{10}})$, and $(0,0,-e^{-\frac{\tau_j}{10}})$. Using Lemma \ref{asymptotics.for.surface} and the convexity of $\bar{M}_{\tau_j}$, we conclude that this triangle lies inside $\bar{M}_{\tau_j}$. Moreover, the angle at $q_j$ converges to $\pi$ as $j \to -\infty$. We now dilate the surface $\bar{M}_{\tau_j}$ to make the mean curvature at $q_j$ equal to $\frac{1}{\sqrt{2}}$. Passing to the limit as $j \to \infty$, we obtain a noncollapsed ancient solution of mean curvature flow which is weakly, but not strictly convex. By the strong maximum principle, the limit splits off a line. By Lemma 3.6 in \cite{Haslhofer-Kleiner1}, the limit is a round cylinder. Therefore, the point $q_j$ lies on an $\varepsilon_0$-neck if $j$ is sufficiently large. This is a contradiction. \\

After these preparations, we now state the main result of this section:

\begin{proposition} 
\label{sharp.asymptotics.for.rescaled.flow}
We have 
\[\sup_{x \in \bar{M}_\tau \cap B_{10}(0)} |x_1^2+x_2^2-2| \leq O(e^{\frac{\tau}{2}}).\] 
\end{proposition} 

\textbf{Proof.} 
In view of Lemma \ref{asymptotics.for.surface}, Lemma \ref{size.of.necklike.region}, and standard interpolation inequalities, we may write $\bar{M}_\tau$ as a graph over the cylinder $\Sigma \cap B_{e^{-\frac{\tau}{100}}}(0)$, and the $C^4$-norm of the height function is bounded by $O(e^{\frac{\tau}{100}})$. We now repeat the argument above, this time with $\rho(\tau) = e^{-\frac{\tau}{1000}}$. As above, we consider the rotated surfaces $\tilde{M}_\tau = S(\tau) \bar{M}_\tau$, where $S(\tau)$ is a function taking values in $SO(3)$. We write each surface $\tilde{M}_\tau$ as a graph over the cylinder, so that 
\[\{x + u(x,\tau) \nu_\Sigma(x): x \in \Sigma \cap B_{2e^{-\frac{\tau}{1000}}}(0)\} \subset \tilde{M}_\tau,\] 
where $\|u(\cdot,\tau)\|_{C^4(\Sigma \cap B_{2e^{-\frac{\tau}{1000}}}(0))} \leq O(e^{\frac{\tau}{200}})$. We choose the matrices $S(\tau)$ in such a way that the orthogonality relations 
\[\int_{\Sigma \cap B_{e^{-\frac{\tau}{1000}}}(0)} e^{-\frac{|x|^2}{4}} \, \langle Ax,\nu_\Sigma \rangle \, u(x,\tau) \, \varphi(e^{\frac{\tau}{1000}} \, x_3) = 0\] 
are satisfied for all $A \in so(3)$. As above, the function $\hat{u}(x,\tau) = u(x,\tau) \, \varphi(e^{\frac{\tau}{1000}} \, x_3)$ satisfies 
\[\Big \| \frac{\partial}{\partial \tau} \hat{u} - \mathcal{L} \hat{u} \Big \|_{\mathcal{H}} \leq O(e^{\frac{\tau}{1000}}) \, \|\hat{u}\|_{\mathcal{H}}.\] 
Hence, if we define 
\begin{align*} 
&U_+(\tau) := \|P_+ \hat{u}(\cdot,\tau)\|_{\mathcal{H}}^2, \\ 
&U_0(\tau) := \|P_0 \hat{u}(\cdot,\tau)\|_{\mathcal{H}}^2, \\ 
&U_-(\tau) := \|P_- \hat{u}(\cdot,\tau)\|_{\mathcal{H}}^2, 
\end{align*} 
then 
\begin{align*} 
&\frac{d}{d\tau} U_+(\tau) \geq U_+(\tau) - O(e^{\frac{\tau}{1000}}) \, (U_+(\tau) + U_0(\tau) + U_-(\tau)), \\ 
&\Big | \frac{d}{d\tau} U_0(\tau) \Big | \leq O(e^{\frac{\tau}{1000}}) \, (U_+(\tau) + U_0(\tau) + U_-(\tau)), \\ 
&\frac{d}{d\tau} U_-(\tau) \leq -U_-(\tau) + O(e^{\frac{\tau}{1000}}) \, (U_+(\tau) + U_0(\tau) + U_-(\tau)). 
\end{align*}
As above, the ODE lemma of Merle and Zaag implies that either $U_0(\tau)+U_-(\tau) \leq o(1) U_+(\tau)$ or $U_+(\tau)+U_-(\tau) \leq o(1) U_0(\tau)$, and the latter case can be ruled out as above using the orthogonality relations and the Brunn-Minkowski inequality. Thus, $U_0(\tau)+U_-(\tau) \leq o(1) U_+(\tau)$. This gives 
\[\frac{d}{d\tau} U_+(\tau) \geq U_+(\tau) - O(e^{\frac{\tau}{1000}}) \, U_+(\tau),\]
hence $U_+(\tau) \leq O(e^\tau)$. This implies $U_0(\tau)+U_-(\tau) \leq o(1) \, U_+(\tau) \leq O(e^\tau)$. From this, we deduce that $\|\hat{u}\|_{\mathcal{H}} \leq O(e^{\frac{\tau}{2}})$. Using Lemma \ref{pde.for.hat.u.2}, we obtain $|A(\tau)| \leq O(e^{\frac{\tau}{2}})$. Since $\lim_{\tau \to -\infty} S(\tau) = \text{\rm id}$, it follows that $|S(\tau)-\text{\rm id}| \leq O(e^{\frac{\tau}{2}})$. Finally, we observe that $u$ satisfies an equation of the form $\frac{\partial}{\partial \tau} u = \tilde{\mathcal{L}} u + \langle A(\tau)x,\nu_\Sigma \rangle$, where $\tilde{\mathcal{L}}$ is an elliptic operator of second order whose coefficients depend on $u$, $\nabla u$, $\nabla^2 u$, and $A(\tau)$. As $\tau \to -\infty$, the coefficients of $\tilde{\mathcal{L}}$ converge smoothly to the corresponding coefficients of $\mathcal{L}$. Using standard interior estimates for parabolic equations, we obtain $\|u(\cdot,\tau)\|_{C^4([0,2\pi] \times [-100,100])} \leq O(e^{\frac{\tau}{2}})$. Since $|S(\tau)-\text{\rm id}| \leq O(e^{\frac{\tau}{2}})$, we conclude that 
\[\sup_{x \in \bar{M}_\tau \cap B_{10}(0)} |x_1^2+x_2^2-2| \leq O(e^{\frac{\tau}{2}}).\] 
This completes the proof of Proposition \ref{sharp.asymptotics.for.rescaled.flow}. \\

\section{Lower bound for $H_{\text{\rm max}}(t)$ as $t \to -\infty$}

\label{barrier}

Let $M_t$, $t \in (-\infty,0]$, be a noncompact ancient solution of mean curvature flow in $\mathbb{R}^3$ which is strictly convex and noncollapsed. For each $t$, we denote by $H_{\text{\rm max}}(t)$ the supremum of the mean curvature of $M_t$.

\begin{proposition}
\label{H_max.finite}
For each $t$, $H_{\text{\rm max}}(t) < \infty$.
\end{proposition}

\textbf{Proof.} 
Let us fix a time $t$ and a small number $\varepsilon$. By Proposition 3.1 in \cite{Haslhofer-Kleiner2}, we can find a compact subset of $M_t$ with the property that every point in the complement of that set lies at the center of an $\varepsilon$-neck. Hence, if $H_{\text{\rm max}}(t) = \infty$, then the surface $M_t$ contains a sequence of $\varepsilon$-necks with radii converging to $0$, but this is impossible in a convex surface. \\

\begin{corollary} 
\label{properties.of.H_max}
The function $H_{\text{\rm max}}(t)$ is continuous and monotone increasing in $t$.
\end{corollary} 

\textbf{Proof.} 
The pointwise curvature derivative estimate of Haslhofer and Kleiner \cite{Haslhofer-Kleiner1},\cite{Haslhofer-Kleiner2} gives $|\frac{\partial}{\partial t} H| \leq CH^3$ for some uniform constant $C$. Consequently, $H_{\text{\rm max}}(t)$ is continuous in $t$. In particular, $H_{\text{\rm max}}(t)$ is uniformly bounded from above on every compact time interval. Hence, we can apply Hamilton's Harnack inequality \cite{Hamilton} to conclude that $H_{\text{\rm max}}(t)$ is monotone increasing in $t$. \\

\begin{proposition}
\label{lower.bound.for.Hmax}
We have $\liminf_{t \to -\infty} H_{\text{\rm max}}(t) > 0$.
\end{proposition} 

\textbf{Proof.} 
Proposition \ref{sharp.asymptotics.for.rescaled.flow} gives 
\[\sup_{x \in (-t)^{-\frac{1}{2}} \, (M_t \cap B_{10(-t)^{\frac{1}{2}}}(0))} |x_1^2+x_2^2-2| \leq O((-t)^{-\frac{1}{2}}).\] 
By assumption, $M_t$ is noncompact and strictly convex. Hence, $M_t$ has exactly one end. Without loss of generality, we may assume that $M_t \cap \{x_3 \leq 0\}$ is compact and $M_t \cap \{x_3 \geq 0\}$ is noncompact. We can find a large constant $K$ such that the curve 
\[(-t)^{-\frac{1}{2}} \, (M_t \cap \{x_3=-2 (-t)^{\frac{1}{2}}\})\] 
lies outside the circle 
\[\{x_1^2+x_2^2=(\sqrt{2}-K (-t)^{-\frac{1}{2}})^2, \, x_3=-2\}\] 
if $-t$ is sufficiently large. Let us consider the self-similar solutions constructed in \cite{Angenent-Daskalopoulos-Sesum}. For $a>0$ large, there exists a surface 
\[\Sigma_a = \{(x_1,x_2,x_3): x_1^2+x_2^2=u_a(-x_3)^2, \, -a \leq x_3 \leq 0\}\] 
which satisfies the shrinker equation $H = \frac{1}{2} \, \langle x,\nu \rangle$. Hence, the surfaces 
\begin{align*} 
\Sigma_{a,t} 
&:= (-t)^{\frac{1}{2}} \, \Sigma_a + (0,0,Ka^2) \\ 
&= \{(x_1,x_2,x_3): x_1^2+x_2^2=(-t) \, u_a((-x_3+Ka^2) (-t)^{-\frac{1}{2}} )^2, \\ 
&\hspace{30mm} Ka^2 - a (-t)^{\frac{1}{2}} \leq x_3 \leq Ka^2\} 
\end{align*} 
evolve by mean curvature flow.

We will use the surfaces $\Sigma_{a,t} \cap \{x_3 \leq -2 (-t)^{\frac{1}{2}}\}$ as barriers for the flow $M_t \cap \{x_3 \leq -2 (-t)^{\frac{1}{2}}\}$. As $t \to -\infty$, the rescaled surfaces $(-t)^{-\frac{1}{2}} \, M_t$ converge in $C_{\text{\rm loc}}^\infty$ to the cylinder $\{x_1^2+x_2^2=2\}$. Moreover, as $t \to -\infty$, the rescaled surfaces $(-t)^{-\frac{1}{2}} \, (\Sigma_{a,t} \cap \{x_3 \leq -2 (-t)^{\frac{1}{2}}\})$ converge to $\Sigma_a \cap \{x_3 \leq -2\}$, which is a compact subset of $\{x_1^2+x_2^2<2\}$. Therefore, $\Sigma_{a,t} \cap \{x_3 \leq -2 (-t)^{\frac{1}{2}}\}$ lies inside $M_t \cap \{x_3 \leq -2 (-t)^{\frac{1}{2}}\}$ if $-t$ is sufficiently large (depending on $a$). 

We next examine the boundary curves $M_t \cap \{x_3 = -2 (-t)^{\frac{1}{2}}\}$ and $\Sigma_{a,t} \cap \{x_3 = -2 (-t)^{\frac{1}{2}}\}$. By our choice of $K$, the curve 
\[(-t)^{-\frac{1}{2}} \, (M_t \cap \{x_3 = -2 (-t)^{\frac{1}{2}}\})\] 
lies outside the circle 
\[\{x_1^2+x_2^2=(\sqrt{2}-K (-t)^{-\frac{1}{2}})^2, \, x_3=-2\}.\] 
Moreover, the curve 
\[(-t)^{-\frac{1}{2}} \, (\Sigma_{a,t} \cap \{x_3 = -2 (-t)^{\frac{1}{2}}\})\] 
is a circle 
\[\{x_1^2+x_2^2=u_a(2 + Ka^2 (-t)^{-\frac{1}{2}})^2, \, x_3=-2\}.\] 
Using Lemma 4.4 in \cite{Angenent-Daskalopoulos-Sesum}, we obtain $u_a(2) \leq \sqrt{2}$ and $u_a(2)-u_a(1) \leq -a^{-2}$ if $a$ is sufficiently large. Moreover, by Lemma 4.2 in \cite{Angenent-Daskalopoulos-Sesum}, the function $u_a: [0,a] \to \mathbb{R}$ is concave. Hence, we obtain 
\begin{align*} 
u_a(2+Ka^2 (-t)^{-\frac{1}{2}}) 
&\leq u_a(2) + Ka^2 (-t)^{-\frac{1}{2}} \, (u_a(2)-u_a(1)) \\ 
&\leq \sqrt{2} - K (-t)^{-\frac{1}{2}} 
\end{align*} 
for $-t \geq 4K^2a^2$. Therefore, the curve $\Sigma_{a,t} \cap \{x_3 = -2 (-t)^{\frac{1}{2}}\}$ lies inside the curve $M_t \cap \{x_3 = -2 (-t)^{\frac{1}{2}}\}$ whenever $-t \geq 4K^2a^2$ and $a$ is sufficiently large. Using the maximum principle, we conclude that the surface $\Sigma_{a,t} \cap \{x_3 \leq -2 (-t)^{\frac{1}{2}}\}$ lies inside the surface $M_t \cap \{x_3 \leq -2 (-t)^{\frac{1}{2}}\}$ whenever $-t \geq 4K^2a^2$ and $a$ is sufficiently large. For $-t=4K^2a^2$, the tip of $\Sigma_{a,t}$ has distance $a (-t)^{\frac{1}{2}} -Ka^2 = Ka^2 = -\frac{t}{4K}$ from the origin. Consequently, the intersection of $M_t$ with the halfline $\{x_1=x_2=0, \, x_3 \leq \frac{t}{4K}\}$ is non-empty if $-t$ is sufficiently large. From this, we deduce that $\limsup_{t \to -\infty} H_{\text{\rm max}}(t) > 0$. Since $H_{\text{\rm max}}(t)$ is monotone increasing in $t$, we conclude that $\liminf_{t \to -\infty} H_{\text{\rm max}}(t) > 0$.

\section{The neck improvement theorem}

\label{neck.improvement}

\begin{definition}
\label{rotation.vf}
Let $K$ be a vector field in $\mathbb{R}^3$. We say that $K$ is a normalized rotation vector field if there exists a matrix $S \in O(3)$ and a point $q \in \mathbb{R}^3$ such that $K(x) = SJS^{-1}(x-q)$, where  
\[J = \begin{bmatrix} 0 & 1 & 0 \\ -1 & 0 & 0 \\ 0 & 0 & 0 \end{bmatrix}.\] 
\end{definition}

Note that we do not require that the origin lies on the axis of rotation.

\begin{lemma}
\label{vector.field.comparison}
There exists a large constant $C$ and small constant $\varepsilon_0>0$ with the following property. Suppose that $M$ is a surface in $\mathbb{R}^3$ and let $\bar{x}$ be a point on $M$. We assume that, after rescaling, the surface $M$ is $\varepsilon_0$-close (in the $C^4$-norm) to a cylinder $S^1 \times [-5,5]$ of radius $1$. Suppose that $K^{(1)}$ and $K^{(2)}$ are normalized rotation vector fields with the following properties: 
\begin{itemize} 
\item $|K^{(1)}| \, H \leq 10$ and $|K^{(2)}| \, H \leq 10$ at the point $\bar{x}$.
\item $|\langle K^{(1)},\nu \rangle| \, H \leq \varepsilon$ and $|\langle K^{(2)},\nu \rangle| \, H \leq \varepsilon$ in a geodesic ball around $\bar{x}$ in $M$ of radius $H(\bar{x})^{-1}$. 
\end{itemize}
Then 
\[\min \bigg \{ \sup_{B_{100 H(\bar{x})^{-1}}(\bar{x})} |K^{(1)}-K^{(2)}|,\sup_{B_{100 H(\bar{x})^{-1}}(\bar{x})} |K^{(1)}+K^{(2)}| \bigg \} \, H(\bar{x}) \leq C\varepsilon.\] 
\end{lemma}

\textbf{Proof.} 
By scaling, we may assume that $H(\bar{x})=1$. We argue by contradiction. If the assertion is false, then there exist a sequence of surfaces $M^{(j)}$, a sequence of points $\bar{x}_j \in M^{(j)}$ satisfying $H(\bar{x}_j)=1$, sequences of normalized rotation vector fields $K^{(1,j)}$ and $K^{(2,j)}$, and a sequences of real number $\varepsilon_j \to 0$ with the following properties: 
\begin{itemize}
\item The surfaces $M^{(j)}$ are $\frac{1}{j}$-close (in the $C^4$-norm) to a cylinder $M = S^1 \times [-5,5]$ of radius $1$. Moreover, we may assume that the axis of the cylinder is the $x_3$-axis. 
\item $|K^{(1,j)}| \leq 10$ and $|K^{(2,j)}| \leq 10$ at the point $\bar{x}_j$. 
\item $|\langle K^{(1,j)},\nu \rangle| \, H \leq \varepsilon_j$ and $|\langle K^{(2,j)},\nu \rangle| \, H \leq \varepsilon_j$ in a geodesic ball around $\bar{x}_j$ in $M$ of radius $1$. 
\item $\sup_{B_{100}(\bar{x}_j)} |K^{(1,j)}-K^{(2,j)}| \geq j \varepsilon_j$. 
\item $\sup_{B_{100}(\bar{x}_j)} |K^{(1,j)}+K^{(2,j)}| \geq j \varepsilon_j$. 
\end{itemize} 
Note that the distance of $\bar{x}_j$ from the axis of rotation of $K^{(1,j)}$ is at most $10$. Hence, after passing to a subsequence if necessary, the vector fields $K^{(1,j)}$ converge to a normalized vector field $K^{(1)}$. Similarly, the vector fields $K^{(2,j)}$ converge to a vector field $K^{(2)}$. Clearly, $K^{(1)}$ and $K^{(2)}$ are tangential to the cylinder $S^1 \times [-5,5]$. Consequently, we have $K^{(1)}(x) = \pm Jx$ and $K^{(2)}(x) = \pm Jx$, where $J$ is defined as above. Without loss of generality, we assume that $K^{(1)}(x)=K^{(2)}(x)=Jx$. For each $j$, we define $\delta_j := \sup_{B_{100}(\bar{x}_j)} |K^{(1,j)}-K^{(2,j)}| \geq j \varepsilon_j$. Clearly, $\delta_j \to 0$. We next consider the Killing vector field $V^{(j)} := \delta_j^{-1} \, (K^{(1,j)}-K^{(2,j)})$. Then $\sup_{B_{100}(\bar{x}_j)} |V^{(j)}| = 1$, and $|\langle V^{(j)},\nu \rangle| \, H \leq 2\delta_j^{-1} \varepsilon_j \leq 2 j^{-1}$ in a geodesic ball around $\bar{x}_j$ in $M$ of radius $1$. Hence, after passing to a subsequence, the vector fields $V^{(j)}$ converge to a non-trivial Killing vector field $V$ on $\mathbb{R}^3$ which is tangential to the cylinder $S^1 \times [-5,5]$. Since $K^{(1,j)}$ and $K^{(2,j)}$ are normalized rotation vector fields, the limit vector field $V$ must be of the form $V(x) = [A,J]x - Jb$ for some matrix $A \in so(3)$ and some vector $b \in \mathbb{R}^3$. However, such a vector field cannot be tangential to the cylinder $S^1 \times [-5,5]$ unless $V$ vanishes identically. This is a contradiction. \\

As in \cite{Huisken-Sinestrari3}, pp.~189--190, we denote by $\mathcal{P}(\bar{x},\bar{t},r,\tau)$ the set of all points $(x,t)$ in space-time such that $x \in B_{g(\bar{t})}(\bar{x},r)$ and $t \in [\bar{t}-\tau,\bar{t}]$. Moreover, we put $\hat{\mathcal{P}}(\bar{x},\bar{t},L,\theta) = \mathcal{P}(\bar{x},\bar{t},L \, H(\bar{x},\bar{t})^{-1},\theta \, H(\bar{x},\bar{t})^{-2})$. We say that $(\bar{x},\bar{t})$ lies on an $\varepsilon$-neck if the parabolic neighborhood $\hat{\mathcal{P}}(\bar{x},\bar{t},100,100)$ is, after rescaling, $\varepsilon$-close (in the $C^{10}$-norm), to a family of shrinking cylinders.

\begin{definition}
Let $M_t$ be a solution of mean curvature flow with positive mean curvature. We say that a point $(\bar{x},\bar{t})$ is $\varepsilon$-symmetric if there exists a normalized rotation vector field $K$ on $\mathbb{R}^3$ such that $|K| \, H \leq 10$ at the point $(\bar{x},\bar{t})$ and $|\langle K,\nu \rangle| \, H \leq \varepsilon$ in the parabolic neighborhood $\hat{\mathcal{P}}(\bar{x},\bar{t},10,100)$.
\end{definition}

Note that the condition that $|K| \, H \leq 10$ at the point $(\bar{x},\bar{t})$ is equivalent to the condition that the distance of the point $\bar{x}$ from the axis of rotation of $K$ is at most $10 \, H(\bar{x},\bar{t})^{-1}$. 

\begin{theorem}[Neck Improvement Theorem]
\label{neck.improvement.theorem}
There exists a large constant $L$ and a small constant $\varepsilon_1$ with the following property. Let $M_t$ be a solution of mean curvature flow, and let $(\bar{x},\bar{t})$ be a point in space-time. Suppose that every point in the parabolic neighborhood $\hat{\mathcal{P}}(\bar{x},\bar{t},L,L^2)$ lies on an $\varepsilon_1$-neck. Moreover, suppose that every point in $\hat{\mathcal{P}}(\bar{x},\bar{t},L,L^2)$ is $\varepsilon$-symmetric, where $\varepsilon \leq \varepsilon_1$. Then $(\bar{x},\bar{t})$ is $\frac{\varepsilon}{2}$-symmetric.
\end{theorem}

\textbf{Proof.} 
Without loss of generality, we assume $\bar{t}=-1$ and $H(\bar{x},-1)=\frac{1}{\sqrt{2}}$. Throughout the proof, we assume that $L$ is sufficiently large, and $\varepsilon_1$ is sufficiently small depending on $L$. In the parabolic neighborhood $\hat{\mathcal{P}}(\bar{x},\bar{t},L,L^2)$, we can approximate $M_t$ by a cylinder $S^1((-2t)^{\frac{1}{2}}) \times \mathbb{R}$, up to errors which are bounded by $C(L)\varepsilon_1$ in the $C^{100}$-norm.

\textit{Step 1:} By assumption, for every point $(x_0,t_0) \in \hat{\mathcal{P}}(\bar{x},-1,L,L^2)$, there exists a normalized vector field $K^{(x_0,t_0)}$ such that $|K^{(x_0,t_0)}| \, H \leq 10$ at the point $(x_0,t_0)$, and $|\langle K^{(x_0,t_0)},\nu \rangle| \, H \leq \varepsilon$ on the parabolic neighborhood $\hat{\mathcal{P}}(x_0,t_0,10,100)$. A repeated application of Lemma \ref{vector.field.comparison} gives 
\[\min \bigg \{ \sup_{B_{10L}(0)} |K^{(\bar{x},-1)} - K^{(x_0,t_0)}|,\sup_{B_{10L}(0)} |K^{(\bar{x},-1)} + K^{(x_0,t_0)}| \bigg \} \leq C(L) \varepsilon\] 
for all $(x_0,t_0) \in \hat{\mathcal{P}}(\bar{x},-1,L,L^2)$. Without loss of generality, we may assume that 
\[\sup_{B_{10L}(0)} |K^{(\bar{x},-1)} - K^{(x_0,t_0)}| \leq C(L)\varepsilon\] 
for all $(x_0,t_0) \in \hat{\mathcal{P}}(\bar{x},-1,L,L^2)$. Moreover, we may assume without loss of generality that $\bar{K} = K^{(\bar{x},-1)}$ is an infinitesimal rotation around the $x_3$-axis, so that $\bar{K}(y) = Jy$, where $J$ is defined as in Definition \ref{rotation.vf}. Finally, we may assume that the point $\bar{x}$ lies in the plane $\{x_3=0\}$.

Let us write $M_t$ as a graph over the $x_3$-axis, so that 
\[\Big \{ (r(\theta,z,t) \cos \theta,r(\theta,z,t) \sin \theta,z): \theta \in [0,2\pi], \, z \in \Big [ -\frac{L}{4},\frac{L}{4} \Big ] \Big \} \subset M_t.\] 
By assumption, the difference $r(\theta,z,t) - (-2t)^{\frac{1}{2}}$ is bounded by $C(L)\varepsilon_1$ in the $C^{100}$-norm. The unit normal vector to $M_t$ is given by 
\begin{align*} 
\nu 
&= \frac{1}{\sqrt{1+r^{-2} \big ( \frac{\partial r}{\partial \theta} \big )^2+\big ( \frac{\partial r}{\partial z} \big )^2}} \\ 
&\cdot \Big [ (\cos \theta,\sin \theta,0) - r^{-1} \, \frac{\partial r}{\partial \theta} \, (-\sin \theta,\cos \theta,0) - \frac{\partial r}{\partial z} \, (0,0,1) \Big ]. 
\end{align*}
We define 
\[u = \langle \bar{K},\nu \rangle = \frac{1}{\sqrt{1+r^{-2} \, \big ( \frac{\partial r}{\partial \theta} \big )^2+\big ( \frac{\partial r}{\partial z} \big )^2}} \, \frac{\partial r}{\partial \theta}.\] 

\textit{Step 2:} For each point $(x_0,t_0) \in \hat{\mathcal{P}}(\bar{x},-1,L,L^2)$, we know that 
\[|\langle K^{(x_0,t_0)},\nu \rangle| \leq C\varepsilon \, (-t_0)^{\frac{1}{2}}\] 
on the parabolic neighborhood $\hat{\mathcal{P}}(x_0,t_0,10,100)$. Moreover, we can find a matrix $S \in O(3)$ and a vector $q \in \mathbb{R}^3$ (depending on $(x_0,t_0)$) such that $\bar{K}(y) - K^{(x_0,t_0)}(y) = Jy-SJS^{-1}(y-q)$ and $|S-\text{\rm id}|+|q| \leq C(L)\varepsilon$. Consequently, there exist real numbers $a_0,a_1,b_0,b_1$ (depending on $(x_0,t_0)$) such that 
\[|a_0|+|a_1|+|b_0|+|b_1| \leq C(L) \varepsilon\] 
and 
\[|\langle \bar{K}-K^{(x_0,t_0)},\nu \rangle -  (a_0+a_1z) \cos \theta - (b_0+b_1z) \sin \theta| \leq C(L)\varepsilon_1\varepsilon\] 
on the parabolic neighborhood $\hat{\mathcal{P}}(x_0,t_0,10,100)$. Note that the numbers $a_0,a_1,b_0,b_1$ account for the fact that the rotation vector fields $K^{(x_0,t_0)}$ and $\bar{K}$ may have different axes of rotation.

Putting these facts together, we obtain 
\[|\langle \bar{K},\nu \rangle - (a_0+a_1z) \cos \theta - (b_0+b_1z) \sin \theta| \leq C\varepsilon (-t_0)^{\frac{1}{2}} + C(L)\varepsilon_1\varepsilon\] 
on the parabolic neighborhood $\hat{\mathcal{P}}(x_0,t_0,10,100)$. To summarize, given any point $(z_0,t_0) \in [-\frac{L}{2},\frac{L}{2}] \times [-\frac{L^2}{4},-1]$, there exist real numbers $a_0,a_1,b_0,b_1$ (depending on $(z_0,t_0)$) such that 
\[|a_0|+|a_1|+|b_0|+|b_1| \leq C(L) \varepsilon\] 
and 
\[|u(\theta,z,t) - (a_0+a_1z) \cos \theta - (b_0+b_1z) \sin \theta| \leq C\varepsilon (-t_0)^{\frac{1}{2}} + C(L)\varepsilon_1\varepsilon\] 
for $z \in [z_0-(-t_0)^{\frac{1}{2}},z_0+(-t_0)^{\frac{1}{2}}]$ and $t \in [2t_0,t_0]$. 

\textit{Step 3:} The function $u = \langle \bar{K},\nu \rangle$ satisfies the evolution equation 
\[\frac{\partial}{\partial t} u = \Delta_{M_t} u + |A|^2 u.\] 
Since $|u| \leq C(L)\varepsilon$, it follows from standard interior estimates for parabolic equations that $|\nabla u| + |\nabla^2 u| \leq C(L) \varepsilon$ for $z \in [-\frac{L}{4},\frac{L}{4}]$ and $t \in [-\frac{L^2}{16},-1]$. Hence, we obtain  
\[\Big | \frac{\partial}{\partial t} u - \frac{\partial^2}{\partial z^2} u - \frac{1}{(-2t)} \, \frac{\partial^2}{\partial \theta^2} u - \frac{1}{(-2t)} \, u \Big | \leq C(L)\varepsilon_1\varepsilon\] 
for $z \in [-\frac{L}{4},\frac{L}{4}]$ and $t \in [-\frac{L^2}{16},-1]$. 

Let $\tilde{u}$ be the solution of the linear equation 
\[\frac{\partial}{\partial t} \tilde{u} = \frac{\partial^2}{\partial z^2} \tilde{u} + \frac{1}{(-2t)} \, \frac{\partial^2}{\partial \theta^2} \tilde{u} + \frac{1}{(-2t)} \, \tilde{u}\] 
in the parabolic cylinder $\{z \in [-\frac{L}{4},\frac{L}{4}], \, t \in [-\frac{L^2}{16},-1]\}$ such that $\tilde{u}=u$ on the parabolic boundary $\{|z|=\frac{L}{4}\} \cup \{t=-\frac{L^2}{16}\}$. Using the maximum principle, we obtain 
\[|u-\tilde{u}| \leq C(L)\varepsilon_1\varepsilon\] 
in the parabolic cylinder $\{z \in [-\frac{L}{4},\frac{L}{4}], \, t \in [-\frac{L^2}{16},-1]\}$. 

\textit{Step 4:} We now analyze the function $\tilde{u}$ using separation of variables. For $m \geq 1$, we put 
\[v_m(z,t) = \frac{1}{\pi} \int_0^{2\pi} \tilde{u}(\theta,z,t) \, \cos(m\theta) \, d\theta\] 
and 
\[w_m(z,t) = \frac{1}{\pi} \int_0^{2\pi} \tilde{u}(\theta,z,t) \, \sin(m\theta) \, d\theta.\] 
The functions $v_m$ and $w_m$ satisfy the evolution equations 
\[\frac{\partial}{\partial t} v_m = \frac{\partial^2}{\partial z^2} v_m + \frac{1-m^2}{(-2t)} \, v_m, \qquad \frac{\partial}{\partial t} w_m = \frac{\partial^2}{\partial z^2} w_m + \frac{1-m^2}{(-2t)} \, w_m.\] 
Consequently, the functions $\hat{v}_m = (-t)^{\frac{1-m^2}{2}} \, v_m$ and $\hat{w}_m = (-t)^{\frac{1-m^2}{2}} \, w_m$ satisfy the linear heat equation 
\[\frac{\partial}{\partial t} \hat{v}_m = \frac{\partial^2}{\partial z^2} \hat{v}_m, \qquad \frac{\partial}{\partial t} \hat{w}_m = \frac{\partial^2}{\partial z^2} \hat{w}_m\] 
for $m \geq 1$. 

\textit{Step 5:} We first consider the modes with $m \geq 2$. For $m \geq 2$, we have 
\[|v_m|+|w_m| \leq (C\varepsilon+C(L)\varepsilon_1\varepsilon) \, (-t)^{\frac{1}{2}}\] 
in the parabolic cylinder $[-\frac{L}{4},\frac{L}{4}] \times [-\frac{L^2}{16},-1]$. This implies 
\[|\hat{v}_m|+|\hat{w}_m| \leq (C\varepsilon+C(L)\varepsilon_1\varepsilon) \, (-t)^{1-\frac{m^2}{2}}\] 
in the parabolic cylinder $[-\frac{L}{4},\frac{L}{4}] \times [-\frac{L^2}{16},-1]$. Using the solution formula for the one-dimensional heat equation with Dirichlet boundary condition on the rectangle $[-\frac{L}{4},\frac{L}{4}] \times [-\frac{L^2}{16},-1]$, we may express $\hat{v}_m(z,t)$ and $\hat{w}_m(z,t)$ as integrals of initial and boundary data. This gives 
\begin{align*} 
|\hat{v}_m(z,t)|+|\hat{w}_m(z,t)| 
&\leq (C \varepsilon+C(L) \varepsilon_1 \varepsilon) \, \Big ( \frac{L}{4} \Big )^{2-m^2}  \\ 
&+ (C\varepsilon+C(L) \varepsilon_1 \varepsilon) L \int_{-\frac{L^2}{16}}^t e^{-\frac{L^2}{100(t-s)}} \, (t-s)^{-\frac{3}{2}} \, (-s)^{1-\frac{m^2}{2}} \, ds \\
&\leq (C \varepsilon+C(L) \varepsilon_1 \varepsilon) \, \Big ( \frac{L}{4} \Big )^{2-m^2} \\ 
&+ (C\varepsilon+C(L) \varepsilon_1 \varepsilon) L^{-1} \int_{-\frac{L^2}{16}}^t (-s)^{\frac{1-m^2}{2}} \, ds \\ 
&\leq (C \varepsilon+C(L) \varepsilon_1 \varepsilon) \, \Big ( \frac{L}{4} \Big )^{2-m^2} \\ 
&+ (C\varepsilon+C(L) \varepsilon_1 \varepsilon) L^{-1} m^{-2} (-t)^{\frac{3-m^2}{2}} 
\end{align*}
for $z \in [-20,20]$ and $t \in [-400,-1]$. Therefore, we obtain  
\begin{align*} 
|v_m(z,t)|+|w_m(z,t)| 
&\leq (C \varepsilon+C(L) \varepsilon_1 \varepsilon) \, \Big ( \frac{L^2}{16(-t)} \Big )^{\frac{2-m^2}{2}} \\ 
&+ (C\varepsilon+C(L) \varepsilon_1 \varepsilon) L^{-1} m^{-2} 
\end{align*}
for $z \in [-20,20]$ and $t \in [-400,-1]$. Summation over $m \geq 2$ yields 
\[\Big | \sum_{m=2}^\infty v_m(z,t) \Big | + \Big | \sum_{m=2}^\infty w_m(z,t) \Big | \leq CL^{-1} \varepsilon + C(L) \varepsilon_1 \varepsilon\] 
for $z \in [-20,20]$ and $t \in [-400,-1]$. 

\textit{Step 6:} We next consider the modes with $m=1$. We have 
\[\frac{\partial}{\partial t} v_1 = \frac{\partial^2}{\partial z^2} v_1, \qquad \frac{\partial}{\partial t} w_1 = \frac{\partial^2}{\partial z^2} w_1.\] 
Moreover, given any point $(z_0,t_0) \in [-\frac{L}{4},\frac{L}{4}] \times [-\frac{L^2}{16},-1]$, there exist constants $a_0,a_1,b_0,b_1$ (depending on $z_0$ and $t_0$) such that 
\[|a_0|,|a_1|,|b_0|,|b_1| \leq C(L)\varepsilon\] 
and 
\[|v_1(z,t)-a_0-a_1z|+|w_1(z,t)-b_0-b_1z| \leq C\varepsilon (-t_0)^{\frac{1}{2}} + C(L)\varepsilon_1 \varepsilon\] 
for $z \in [z_0-(-t_0)^{\frac{1}{2}},z_0+(-t_0)^{\frac{1}{2}}]$ and $t \in [2t_0,t_0]$. Using standard interior estimates for the linear heat equation, we obtain 
\[\Big | \frac{\partial^2 v_1}{\partial z^2} \Big | + \Big | \frac{\partial^2 w_1}{\partial z^2} \Big | \leq (C\varepsilon + C(L)\varepsilon_1\varepsilon) \, (-t)^{-\frac{1}{2}}\] 
in the parabolic cylinder $[-\frac{L}{4},\frac{L}{4}] \times [-\frac{L^2}{16},-1]$. As above, we can use the solution formula for the one-dimensional heat equation with Dirichlet boundary condition on the rectangle $[-\frac{L}{4},\frac{L}{4}] \times [-\frac{L^2}{16},-1]$ to express $\frac{\partial^2 v_1}{\partial z^2}$ and $\frac{\partial^2 w_1}{\partial z^2}$ as integrals of initial and boundary data. This gives 
\begin{align*} 
&\Big | \frac{\partial^2 v_1}{\partial z^2}(z,t) \Big |+\Big | \frac{\partial^2 w_1}{\partial z^2}(z,t) \Big | \\ 
&\leq (C\varepsilon + C(L)\varepsilon_1\varepsilon) \, \Big ( \frac{L}{4} \Big )^{-1} \\ 
&+ (C\varepsilon+C(L) \varepsilon_1 \varepsilon) L \int_{-\frac{L^2}{16}}^t e^{-\frac{L^2}{100(t-s)}} \, (t-s)^{-\frac{3}{2}} \, (-s)^{-\frac{1}{2}} \, ds \\ 
&\leq (C\varepsilon + C(L)\varepsilon_1\varepsilon) \, \Big ( \frac{L}{4} \Big )^{-1} \\ 
&+ (C\varepsilon+C(L) \varepsilon_1 \varepsilon) L^{-2} \int_{-\frac{L^2}{16}}^t (-s)^{-\frac{1}{2}} \, ds \\ 
&\leq (C\varepsilon + C(L)\varepsilon_1\varepsilon) \, L^{-1} 
\end{align*} 
for $z \in [-20,20]$ and $t \in [-400,-1]$. Consequently, we can find real numbers $A_0,A_1,B_0,B_1$ such that 
\[|v_1(z,t) - A_0-A_1z|+|w_1(z,t) - B_0-B_1z| \leq CL^{-1} \varepsilon + C(L) \varepsilon_1\varepsilon\] 
for $z \in [-20,20]$ and $t \in [-400,-1]$. 

\textit{Step 7:} Finally, we consider the mode with $m=0$. Using identity 
\[\int_0^{2\pi} u \, \sqrt{1+r^{-2} \, \big ( \frac{\partial r}{\partial \theta} \big )^2+\big ( \frac{\partial r}{\partial z} \big )^2} \, d\theta = \int_0^{2\pi} \frac{\partial r}{\partial \theta} \, d\theta = 0\] 
together with the estimates $|u| \leq C(L)\varepsilon$, $|\frac{\partial r}{\partial \theta}|+|\frac{\partial r}{\partial z}| \leq C(L)\varepsilon_1$, we obtain 
\[\Big | \int_0^{2\pi} u(\theta,z,t) \, d\theta \Big | \leq C(L) \varepsilon_1 \varepsilon,\] 
hence 
\[\Big | \int_0^{2\pi} \tilde{u}(\theta,z,t) \, d\theta \Big | \leq C(L) \varepsilon_1 \varepsilon\] 
for $z \in [-20,20]$ and $t \in [-400,-1]$. 

\textit{Step 8:} To summarize, we have shown that there exist $A_0,A_1,B_0,B_1$ such that 
\[|\tilde{u} - (A_0+A_1z) \cos \theta - (B_0+B_1z) \sin \theta| \leq CL^{-1} \varepsilon + C(L) \varepsilon_1 \varepsilon\] 
for $z \in [-20,20]$ and $t \in [-400,-1]$. This directly implies 
\[|u - (A_0+A_1z) \cos \theta - (B_0+B_1z) \sin \theta| \leq CL^{-1} \varepsilon + C(L) \varepsilon_1 \varepsilon\] 
for $z \in [-20,20]$ and $t \in [-400,-1]$. In particular, $|A_0|+|A_1|+|B_0|+|B_1| \leq C(L)\varepsilon$. Hence, there exists a normalized rotation vector field $\tilde{K}$ such that 
\[|\langle \tilde{K},\nu \rangle| \leq CL^{-1} \varepsilon + C(L) \varepsilon_1 \varepsilon\] 
in the parabolic neighborhood $\hat{\mathcal{P}}(\bar{x},-1,10,100)$. Therefore, $(\bar{x},-1)$ is $(CL^{-1} \varepsilon + C(L) \varepsilon_1 \varepsilon)$-symmetric. In particular, if we choose $L$ sufficiently large and $\varepsilon_1$ sufficiently small (depending on $L$), then $(\bar{x},-1)$ is $\frac{\varepsilon}{2}$-symmetric.

\section{Proof of rotational symmetry}

\label{rotational.symmetry}

Let $M_t$, $t \in (-\infty,0]$, be a noncompact ancient solution of mean curvature flow in $\mathbb{R}^3$ which is strictly convex and noncollapsed.

\begin{lemma}
\label{point.of.maximum.curvature}
If $-t$ is sufficiently large, then there exists a unique point $p_t \in M_t$, where the mean curvature attains its maximum. Moreover, the Hessian of the mean curvature at $p_t$ is negative definite. In particular, $p_t$ varies smoothly in $t$. 
\end{lemma} 

\textbf{Proof.}
We know that $M_t \cap B_{8(-t)^{\frac{1}{2}}}(0)$ is a neck with radius $(-2t)^{\frac{1}{2}}$. The complement $M_t \setminus B_{8(-t)^{\frac{1}{2}}}(0)$ has two connected components, one of which is compact and one of which is noncompact. On the noncompact connected component, the mean curvature is bounded from above by $C \, (-t)^{-\frac{1}{2}}$. On the other hand,  we have shown in Proposition \ref{lower.bound.for.Hmax} that $H_{\text{\rm max}}(t)$ is bounded away from $0$. Consequently, if $-t$ is sufficiently large, then the maximum of the mean curvature is attained at some point $p_t \in M_t$.

We next consider an arbitrary sequence $t_j \to -\infty$, and define $M_t^{(j)} := M_{t+t_j}-p_{t_j}$, where $p_{t_j}$ is the point on $M_{t_j}$ where the mean curvature attains its maximum. After passing to a subsequence if necessary, the sequence $M_t^{(j)}$ converges in $C_{\text{\rm loc}}^\infty$ to a smooth eternal solution. Moreover, there exists a point on the limit solution where the mean curvature attains its space-time maximum. By work of Hamilton \cite{Hamilton}, the limit solution must be a translating soliton. By \cite{Haslhofer}, the limit is the bowl soliton. Hence, if $-t$ is sufficiently large, then $p_t$ is the only point on $M_t$ where the maximum of the mean curvature is attained, and the Hessian of the mean curvature at $p_t$ is negative definite. This completes the proof of Lemma \ref{point.of.maximum.curvature}. \\

Let $\varepsilon_1$ and $L$ be the constants in the Neck Improvement Theorem. Since $H_{\text{\rm max}}(t)$ is uniformly bounded from below, Proposition 3.1 in \cite{Haslhofer-Kleiner2} implies that there exists a large constant $\Lambda$ with the property that every point $x \in M_t$ with $|x-p_t| \geq \Lambda$ lies at the center of an $\varepsilon_1$-neck and satisfies $H(x,t) \, |x-p_t| \geq 10^6 \, L$. 

\begin{lemma}
\label{distance.to.tip.decreases}
There exists a time $T<0$ with the following property: Suppose that $\bar{t} \leq T$, and $\bar{x}$ is a point on $M_{\bar{t}}$ satisfying $|\bar{x}-p_{\bar{t}}| \geq \Lambda$. Then $|\bar{x}-p_t| \geq |\bar{x}-p_{\bar{t}}|$ for all $t \leq \bar{t}$.
\end{lemma} 

\textbf{Proof.} 
If $-t$ is sufficiently large, then $M_t$ looks like the bowl soliton near the point $p_t$. Hence, if $-t$ is sufficiently large, then the vector $\frac{d}{dt} p_t$ is almost parallel to $-\nu(p_t,t)$. Consequently, we can find a time $T<0$ with the property that $\langle x-p_t,\frac{d}{dt} p_t \rangle > 0$ whenever $t \leq T$ and $|x-p_t| \geq \Lambda$. This implies $\frac{d}{dt} |x-p_t| = -\langle \frac{x-p_t}{|x-p_t|},\frac{d}{dt} p_t \rangle < 0$ whenever $t \leq T$ and $|x-p_t| \geq \Lambda$.

We will show that $T$ has the desired property. To prove this, we consider a time $\bar{t} \leq T$ and a point $\bar{x} \in M_{\bar{t}}$ such that $|\bar{x}-p_{\bar{t}}| \geq \Lambda$. We claim that $|\bar{x}-p_t| \geq |\bar{x}-p_{\bar{t}}|$ for all $t \leq \bar{t}$. Indeed, if this is false, then we define $\tilde{t} := \sup \{t \leq \bar{t}: |\bar{x}-p_t| < |\bar{x}-p_{\bar{t}}|\}$. Clearly, $\tilde{t} < \bar{t}$, and $|\bar{x}-p_t| \geq |\bar{x}-p_{\bar{t}}| \geq \Lambda$ for all $t \in [\tilde{t},\bar{t}]$. In view of our choice of $T$, we obtain $\frac{d}{dt} |\bar{x}-p_t| < 0$ for all $t \in [\tilde{t},\bar{t}]$. Consequently, $|\bar{x}-p_{\tilde{t}}| > |\bar{x}-p_{\bar{t}}|$, which contradicts the definition of $\tilde{t}$. This completes the proof of Lemma \ref{distance.to.tip.decreases}. \\

\begin{proposition} 
\label{iteration}
If $t \leq T$, $x \in M_t$ and $|x-p_t| \geq 2^{\frac{j}{400}} \, \Lambda$, then $(x,t)$ is $2^{-j} \varepsilon_1$-symmetric. 
\end{proposition}

\textbf{Proof.} 
We argue by induction on $j$. For $j=0$, the assertion is true. Suppose now that $j \geq 1$ and the assertion holds for $j-1$. We claim that the assertion holds for $j$. Suppose this is false. Then there exists a time $\bar{t} \leq T$ and a point $\bar{x} \in M_{\bar{t}}$ such that $|\bar{x}-p_{\bar{t}}| \geq 2^{\frac{j}{400}} \, \Lambda$ and $(\bar{x},\bar{t})$ is not $2^{-j} \varepsilon_1$-symmetric. By the Neck Improvement Theorem, there exists a point $(x,t) \in \hat{\mathcal{P}}(\bar{x},\bar{t},L,L^2)$ such that either $(x,t)$ is not $2^{-j+1} \varepsilon_1$-symmetric or $(x,t)$ does not lie at the center of an $\varepsilon_1$-neck. In view of the induction hypothesis, we conclude that $|x-p_t| \leq 2^{\frac{j-1}{400}} \, \Lambda$. Since $t \leq \bar{t} \leq T$, Lemma \ref{distance.to.tip.decreases} gives $|\bar{x}-p_{\bar{t}}| \leq |\bar{x}-p_t|$. Putting these facts together, we obtain 
\begin{align*} 
|\bar{x}-p_{\bar{t}}| 
&\leq |\bar{x}-p_t| \\ 
&\leq |x-p_t|+|x-\bar{x}| \\ 
&\leq 2^{\frac{j-1}{400}} \, \Lambda + 10 \, L \, H(\bar{x},\bar{t})^{-1} \\ 
&\leq 2^{-\frac{1}{400}} \, |\bar{x}-p_{\bar{t}}| + \frac{1}{1000} \, |\bar{x}-p_{\bar{t}}| \\ 
&< |\bar{x}-p_{\bar{t}}|. 
\end{align*}
This is a contradiction. \\

\begin{theorem} 
\label{symmetry}
The surface $M_t$ is rotationally symmetric for each $t \leq T$.
\end{theorem}

\textbf{Proof.} 
Let us fix a time $\bar{t} \leq T$. For $j$ sufficiently large, we denote by $\Omega^{(j)}$ the set of all points $(x,t)$ in space-time such that  $t \leq \bar{t}$ and $|x-p_t| \leq 2^{\frac{j}{400}} \, \Lambda$. If $j$ is sufficiently large, then $H(x,t) \geq 2 \cdot 2^{-\frac{j}{400}}$ for each point $(x,t) \in \Omega^{(j)}$. By Proposition \ref{iteration}, every point $(x,t) \in \partial \Omega^{(j)}$ is $2^{-j} \varepsilon_1$-symmetric. In other words, for each point $(x,t) \in \partial \Omega^{(j)}$, there exists a normalized rotation vector field $K^{(x,t)}$ such that $|\langle K^{(x,t)},\nu \rangle| \, H \leq 2^{-j} \varepsilon_1$ on $\hat{\mathcal{P}}(x,t,10,100)$. Using Lemma \ref{vector.field.comparison}, we can control how the axis of rotation of $K^{(x,t)}$ varies as we vary the point $(x,t)$. More precisely, if $(x_1,t_1)$ and $(x_2,t_2)$ are two points on $\partial \Omega^{(j)}$ satisfying $(x_2,t_2) \in \hat{\mathcal{P}}(x_1,t_1,1,1)$, then 
\begin{align*} 
&\min \bigg \{ \sup_{B_{10 H(x_2,t_2)^{-1}}(x_2)} |K^{(x_1,t_1)} - K^{(x_2,t_2)}|,\sup_{B_{10 H(x_2,t_2)^{-1}}(x_2)} |K^{(x_1,t_1)} + K^{(x_2,t_2)}| \bigg \} \\ 
&\leq C \, 2^{-j} \, H(x_2,t_2)^{-1}. 
\end{align*}
Hence, there exists a single normalized rotation vector field $K^{(j)}$ with the following property: if $(x,t)$ is a point in $\partial \Omega^{(j)}$ satisfying $\bar{t}-2^{\frac{j}{100}} \leq t \leq \bar{t}$, then 
\[\min \{ |K^{(x,t)} - K^{(j)}|,|K^{(x,t)} + K^{(j)}|\} \leq C \, 2^{-\frac{j}{2}}\] 
at the point $(x,t)$. This implies $|\langle K^{(j)},\nu \rangle| \leq C \, 2^{-\frac{j}{2}}$ for all points $(x,t) \in \partial \Omega^{(j)}$ satisfying $\bar{t}-2^{\frac{j}{100}} \leq t \leq \bar{t}$. Moreover, we clearly have $|\langle K^{(j)},\nu \rangle| \leq C \, 2^{\frac{j}{100}}$ for all points $(x,t) \in \Omega^{(j)}$ with $t=\bar{t}-2^{\frac{j}{100}}$. 

We now define a function $f^{(j)}: \Omega^{(j)} \to \mathbb{R}$ by 
\[f^{(j)} := \exp(-2^{-\frac{j}{200}} (\bar{t}-t)) \, \frac{\langle K^{(j)},\nu \rangle}{H - 2^{-\frac{j}{400}}}.\] 
Using the estimate for $\langle K^{(j)},\nu \rangle$, we obtain  
\[|f^{(j)}(x,t)| \leq \frac{C\, 2^{-\frac{j}{2}}}{2\cdot 2^{-\frac{j}{400}}-2^{-\frac{j}{400}}} \leq C\,2^{-\frac{j}{4}}\]
for all points $(x,t) \in \partial \Omega^{(j)}$ satisfying $\bar{t}-2^{\frac{j}{100}} \leq t \leq \bar{t}$, and
\[|f^{(j)}(x,t)| \leq \exp(-2^{-\frac{j}{200}+\frac{j}{100}} ) \cdot \frac{C\, 2^{\frac{j}{100}}}{2\cdot 2^{-\frac{j}{400}}-2^{-\frac{j}{400}}} \leq C\,2^{-\frac{j}{4}}\]
for all points $(x,t) \in \Omega^{(j)}$ with $t=\bar{t}-2^{\frac{j}{100}}$. Using the identities 
\[\frac{\partial}{\partial t} H = \Delta H + |A|^2 H\] 
and 
\[\frac{\partial}{\partial t} \langle K^{(j)},\nu \rangle = \Delta \langle K^{(j)},\nu \rangle + |A|^2 \langle K^{(j)},\nu \rangle,\] 
we compute 
\[\frac{\partial}{\partial t} f^{(j)} = \Delta f^{(j)} + \frac{2}{H-2^{-\frac{j}{400}}} \, \langle \nabla H,\nabla f^{(j)} \rangle - 2^{-\frac{j}{400}} \, \Big ( \frac{|A|^2}{H-2^{-\frac{j}{400}}} - 2^{-\frac{j}{400}} \Big ) \, f^{(j)}.\] 
On the set $\Omega^{(j)}$, we have 
\[\frac{|A|^2}{H-2^{-\frac{j}{400}}} - 2^{-\frac{j}{400}} \geq \frac{1}{2} \, \frac{H^2}{H-2^{-\frac{j}{400}}} - 2^{-\frac{j}{400}} \geq \frac{1}{2} \, H - 2^{-\frac{j}{400}} \geq 0.\] 
Using the maximum principle, we conclude that 
\begin{align*} 
&\sup_{(x,t) \in \Omega^{(j)}, \, \bar{t}-2^{\frac{j}{100}} \leq t \leq \bar{t}} |f^{(j)}(x,t)| \\ 
&\leq \max \bigg \{ \sup_{(x,t) \in \partial \Omega^{(j)}, \, \bar{t}-2^{\frac{j}{100}} \leq t \leq \bar{t}} |f^{(j)}(x,t)|,\sup_{(x,t) \in \Omega^{(j)}, \, t=\bar{t}-2^{\frac{j}{100}}} |f^{(j)}(x,t)| \bigg \} \\ 
&\leq C \, 2^{-\frac{j}{4}}. 
\end{align*} 
This gives $|\langle K^{(j)},\nu \rangle| \leq C \, 2^{-\frac{j}{4}}$ for all points $(x,t) \in \Omega^{(j)}$ with $t=\bar{t}$. From this, we deduce that the distance of the point $p_{\bar{t}}$ from the axis of rotation of $K^{(j)}$ is bounded from above by a uniform constant which is independent of $j$. Hence, if we send $j \to \infty$, the vector fields $K^{(j)}$ converge to a normalized rotation vector field in $\mathbb{R}^3$ which is tangential along $M_{\bar{t}}$. This completes the proof of Theorem \ref{symmetry}. \\

Once we know that $M_t$ is rotationally symmetric for $-t$ sufficiently large, it follows from standard arguments that $M_t$ is rotationally symmetric for all $t$:

\begin{proposition}
\label{symmetry.preserved}
Suppose that $M_{\bar{t}}$ is rotationally symmetric for some $\bar{t} \in (-\infty,0]$. Then, for each $t \in [\bar{t},0]$, $M_t$ is rotationally symmetric with the same axis.
\end{proposition}

\textbf{Proof.} 
By Proposition \ref{H_max.finite} and Corollary \ref{properties.of.H_max}, the flow $M_t$, $t \in (-\infty,0]$, has bounded curvature. Without loss of generality, we may assume that $\sup_{M_t} |A|^2 \leq 2$ for each $t \in (-\infty,0]$. If $K$ is a rotation vector field in $\mathbb{R}^3$, then 
\[\frac{\partial}{\partial t} \langle K,\nu \rangle = \Delta \langle K,\nu \rangle + |A|^2 \langle K,\nu \rangle.\] 
Moreover, since $|A|^2 \leq 2$, the function $\rho(x,t) := e^{8t} (|x|^2+1)$ satisfies 
\[\frac{\partial}{\partial t} \rho > \Delta \rho + |A|^2 \rho\] 
for $t \in (-\infty,0]$. By the maximum principle, the quantity $\sup_{M_t} \frac{|\langle K,\nu \rangle|}{\rho}$ is monotone decreasing for $t \in (-\infty,0]$. In particular, if $\langle K,\nu \rangle = 0$ at each point on $M_{\bar{t}}$, then $\langle K,\nu \rangle = 0$ on $M_t$ for all $t \in [\bar{t},0]$. \\

\section{Uniqueness of ancient solutions with rotational symmetry}

\label{analysis.in.rotationally.symmetric.case}

Let $M_t$ be an ancient solution satisfying the assumptions of Theorem \ref{main.thm}. By the results in Section \ref{rotational.symmetry}, $M_t$ is rotationally symmetric. Without loss of generality, we may assume that $M_t$ is symmetric with respect to the $x_3$-axis. Thus, there exists a function $f(r,t)$ such that the solution $M_t$ consists of the points $(r\cos \theta,r\sin\theta,f(r,t)) \in \mathbb{R}^3$. Moreover, the function $f(r,t)$ satisfies the following evolution equation:
\[f_t= \frac{f_{rr}}{1+f_r^2}+\frac{1}{r}f_r.\]
Note that $f(r,t)$ may not be defined for all $r$.

Next, we can consider the radius $r$ as a function of $(z,t)$. Namely, the radius function $r(z,t)$ is defined by 
\[f \big(r(z,t),t \big )=z.\]
Then $r(z,t)$ satisfies the following equation (see also \cite{Angenent-Daskalopoulos-Sesum}):
\[r_t =\frac{r_{zz}}{1+r_z^2} -\frac{1}{r}.\]
Note that the convexity of $M_t$ yields
\begin{align*} 
&r>0,&&r_z>0, &&r_t < 0 ,&& r_{zz} <0. 
\end{align*}
Without loss of generality, we assume that the tip of $M_0$ is at the origin. In other words, $f(0,0)=0$ and $r(0,0)=0$.

Let $q_t = (0,0,f(0,t))$ denote the tip of $M_t$, and let $H_{\text{\rm tip}}(t)$ denote the mean curvature of $M_t$ at the tip $q_t$. Using the Harnack inequality, we conclude that $H_{\text{\rm tip}}(t)$ is monotone increasing. In particular, the limit $\mathcal{H} := \lim_{t \to -\infty} H_{\text{\rm tip}}(t)$ exists. Using Proposition \ref{lower.bound.for.Hmax}, we obtain $\mathcal{H}>0$.

We first prove that $f_t(r,t)$ is monotone increasing in $t$.

\begin{proposition}
\label{consequence.of.harnack}
We have $f_{tt}(r,t)\geq 0$ everywhere.
\end{proposition}

\textbf{Proof.} 
Hamilton's Harnack inequality \cite{Hamilton} implies that 
\[\frac{\partial}{\partial t} H+2 V^i \nabla_i H +h_{ij}V^iV^j \geq 0\] 
for every vector field $V$.

Let $\omega = (0,0,-1)$ denote the vertical vector field in $\mathbb{R}^3$, and let $V = -H \, \langle \omega,\nu \rangle^{-1} \, \omega^{\text{\rm tan}}$. For this choice of $V$, the Harnack inequality takes the form 
\[\Big ( \frac{\partial}{\partial t} + V^i \, \nabla_i \Big ) (H \, \langle \omega,\nu \rangle^{-1}) \geq 0.\] 
A straightforward calculation gives 
\[f_t(r,t) = H \, \langle \omega,\nu \rangle^{-1}\] 
and 
\[f_{tt}(r,t) = \Big ( \frac{\partial}{\partial t} + V^i \, \nabla_i \Big ) (H \, \langle \omega,\nu \rangle^{-1}).\] 
Putting these facts together, the assertion follows. \\

We next show that $f_t(r,t)$ is bounded from below.

\begin{proposition}
\label{lower.bound.for.f_t}
We have $f_t(r,t) \geq \mathcal{H}$ at each point in space-time. Moreover, for each $r_0 > 0$,
\[\lim_{t \to -\infty} \sup_{r \leq r_0} f_t(r,t) = \mathcal{H}.\]
\end{proposition}

\textbf{Proof.}
We consider an arbitrary sequence $t_j \to -\infty$, and define $M_t^{(j)} := M_{t+t_j}-q_{t_j}$. We apply the compactness theorem for ancient solutions (cf. \cite{Haslhofer-Kleiner1}, Theorem 1.10) to the sequence $M_t^{(j)}$. Hence, after passing to a subsequence if necessary, the sequence $M_t^{(j)}$ converges in $C_{\text{\rm loc}}^\infty$ to a smooth eternal solution, which is rotationally symmetry. Moreover, on the limit solution, the mean curvature at the tip equals $\mathcal{H}$ at all times. Hence, equality holds in the Harnack inequality. By work of Hamilton \cite{Hamilton}, the limit solution must be a soliton which is translating with speed $\mathcal{H}$. This directly implies 
\[\lim_{j \to \infty} \sup_{r \leq r_0} |f_t(r,t_j) - \mathcal{H}| = 0\]
for every $r_0>0$. Since $f_{tt}(r,t) \geq 0$ by Proposition \ref{consequence.of.harnack}, we conclude that $f_t(r,t) \geq \mathcal{H}$ for all $r$ and $t$. \\

We next prove that $f_t(r,t)$ is monotone increasing in $r$.

\begin{proposition}
\label{monotonicity.in.r}
We have $f_{tr}(r,t) \geq 0 $ everywhere. 
\end{proposition}

\textbf{Proof.} 
Let us fix a time $t_0$ and a radius $r_0$ such that $f(r_0,t_0)$ is defined. For each $T < t_0$, we consider the parabolic region $Q_T=\{x_1^2+x_2^2 \leq r_0^2, t \in [T,t_0]\}$. Using the equations 
\[\frac{\partial}{\partial t} H = \Delta H + |A|^2 H\] 
and 
\[\frac{\partial}{\partial t} \langle \omega,\nu \rangle = \Delta \langle \omega,\nu \rangle + |A|^2 \langle \omega,\nu \rangle,\] 
we conclude that the maximum $\sup_{Q_T} H \, \langle \omega,\nu \rangle^{-1}$ must be attained on the parabolic boundary of $Q_T$. This gives 
\begin{align*} 
&\sup_{x_1^2+x_2^2 \leq r_0^2, t=t_0} H \, \langle \omega,\nu \rangle^{-1} \\ 
&\leq \max \Big \{\sup_{x_1^2+x_2^2 = r_0^2, T \leq t \leq t_0} H \, \langle \omega,\nu \rangle^{-1}, \, \sup_{x_1^2+x_2^2 \leq r_0^2, t=T} H \, \langle \omega,\nu \rangle^{-1} \Big \}. 
\end{align*}
Since $f_t(r,t) = H \, \langle \omega,\nu \rangle^{-1}$, it follows that 
\begin{align*} 
\sup_{r \leq r_0} f_t(r,t_0) 
&\leq \max \Big \{ \sup_{T \leq t \leq t_0} f_t(r_0,t), \, \sup_{r \leq r_0} f_t(r,T) \Big \} \\ 
&= \max \Big \{ f_t(r_0,t_0), \, \sup_{r \leq r_0} f_t(r,T) \Big \}, 
\end{align*}
where in the last step we have used Proposition \ref{consequence.of.harnack}. We now send $T \to -\infty$. Since $\lim_{t \to -\infty} \sup_{r \leq r_0} f_t(r,t) = \mathcal{H}$, we conclude that 
\[\sup_{r \leq r_0} f_t(r,t_0) \leq \max \{f_t(r_0,t_0),\mathcal{H}\} = f_t(r_0,t_0).\] 
This completes the proof of Proposition \ref{monotonicity.in.r}. \\

We recall that $M_t$ is strictly convex and noncollapsed and $H_{\text{\rm tip}}(t)$ is bounded from below by $\mathcal{H}$. By Proposition 3.1 in \cite{Haslhofer-Kleiner2}, there exists a small constant $\varepsilon_0 \in (0,\frac{1}{20})$ and a decreasing function $\Lambda: (0,\varepsilon_0] \to \mathbb{R}$ such that given any $\varepsilon \in (0,\varepsilon_0]$, if $|\bar{x}-q_t| \geq \Lambda(\varepsilon)$, then $(\bar{x}, \bar{t})$ is a center of $\varepsilon$-neck. (Alternatively, this can be deduced from Theorem 7.14 and Lemma 7.4 in \cite{Huisken-Sinestrari3}.) 

\begin{lemma}
\label{bound.for.rr_z}
On every $\varepsilon_0$-neck, $rr_z=\frac{r}{f_r} \leq (1+2\varepsilon_0)\mathcal{H}^{-1}$. 
\end{lemma}

\textbf{Proof.} 
On an $\varepsilon_0$-neck, we have we have $\frac{1}{f_r} = r_z \leq \varepsilon_0$. Moreover, the principal curvature in radial direction is bounded by $\frac{\varepsilon_0}{r}$. This gives 
\[\frac{f_{rr}}{(1+f_r^2)^{\frac{3}{2}}} \leq \frac{\varepsilon_0}{r}.\] 
Using Proposition \ref{lower.bound.for.f_t}, we obtain 
\[\mathcal{H} \leq f_t = \frac{f_{rr}}{1+f_r^2}+\frac{f_r}{r}\leq \frac{\varepsilon_0}{r} \, (1+f_r^2)^{\frac{1}{2}}+\frac{f_r}{r} \leq (1+2\varepsilon_0)\frac{f_r}{r},\] 
as claimed. \\

\begin{lemma}
\label{interior.estimate.for.r}
There exists a constant $C_0 \geq 1$ such that $r^m \, \big | \frac{\partial^m}{\partial z^m} r \big | \leq C_0$ holds for $m=1,2,3$ at center of $\varepsilon_0$-necks with $r \geq 1$. 
\end{lemma}

\textbf{Proof.}
For $m=1$, the assertion follows from Lemma \ref{bound.for.rr_z}. Let $u = \langle \omega,\nu \rangle$, where $\omega = (0,0,-1)$. Then 
\[\frac{\partial}{\partial t} u = \Delta u  + |A|^2 u.\] 
Moreover, $u \leq \mathcal{H}^{-1} \, H$ by Proposition \ref{lower.bound.for.f_t}. Standard interior estimates imply that 
\[|\nabla^m u|^2\leq C \, H^{2m+2}\]  
for $m=1,2$ at the center of an $\varepsilon_0$-necks.
 
In the parametrization $(z,\theta) \mapsto (r(z) \cos\theta,r(z)\sin\theta,z)$, the induced metric is given by $g_{zz}=1+r_z^2$, $g_{z\theta}=0$, $g_{\theta\theta}=r^2$. Moreover, $u = r_z \, (1+r_z^2)^{-\frac{1}{2}}$ and $u_z=r_{zz}(1+r_z^2)^{-\frac{3}{2}}$. Hence, $|\nabla u|^2=g^{zz} u_z^2=r_{zz}^2(1+r_z^2)^{-4}$. In addition, $r_z \leq \varepsilon_0$ and $Hr \leq 1+\varepsilon_0$ hold in every $\varepsilon_0$-neck. Therefore, the inequality $|\nabla u|^2 \leq C \, H^4 \leq C \, r^{-4}$ implies $r^4 r_{zz}^2 \leq C$. Similarly, $|\nabla^2 u|^2 \leq C \, H^6 \leq C \, r^{-6}$ gives $r^6 r_{zzz}^2 \leq C$. \\

\begin{proposition}
\label{estimate.for.r_zz}
Let $C_1=2+2\Lambda(\varepsilon_0)+9\mathcal{H}^{-2}$. If $r \geq C_1$, then $0 \leq -r_{zz}(z,t) \leq C_2 r(z,t)^{-\frac{5}{2}}$.
\end{proposition}

\textbf{Proof.} 
Clearly, $-r_{zz} \geq 0$ since $M_t$ is convex. To prove the upper bound for $-r_{zz}$, let us fix a point $(\bar{r}, \bar{t})$ satisfying $\bar{r} \geq C_1 \geq 2$, and let $\bar{z}=f(\bar{r},\bar{t})$. Then we have  $\frac{1}{2}\bar{r} \geq \frac{1}{2}C_1 \geq \Lambda(\varepsilon_0)$ by definition of $C_1$. Hence, every point $(x,t)$ with $r=(x_1^2+x_2^2)^\frac{1}{2} \geq \frac{1}{2}\bar{r}$ lies at the center of an $\varepsilon_0$-neck. 

Using Lemma \ref{bound.for.rr_z} and $\varepsilon_0 \leq \frac{1}{20}$, we obtain
\begin{align*}
\bar{z}-f \Big ( \frac{\bar{r}}{2},\bar{t} \Big ) =\int_{\frac{\bar{r}}{2}}^{\bar{r}} f_r(r,\bar{t}) \, dr \geq \int_{\frac{\bar{r}}{2}}^{\bar{r}} \frac{\mathcal{H}  }{1+2\varepsilon_0} \, r \, dr \geq \frac{1}{3}\mathcal{H} \bar{r}^2.
\end{align*}
Since $\bar{r} \geq C_1 \geq 9\mathcal{H}^{-2}$, it follows that $f(\frac{\bar{r}}{2},\bar{t}) \leq \bar{z} -\bar{r}^{\frac{3}{2}}$. In other words, $r(z,\bar{t}) \geq \frac{\bar{r}}{2}$ for $z \in [\bar{z}-\bar{r}^{\frac{3}{2}},\bar{z}+\bar{r}^{\frac{3}{2}}]$. Since $r(z,t)$ is decreasing in $t$, it follows that $r(z,t) \geq \frac{\bar{r}}{2}$ for $(z,t) \in Q = [\bar{z}-\bar{r}^{\frac{3}{2}},\bar{z}+\bar{r}^{\frac{3}{2}}] \times [\bar{t}-\bar{r}^{3},\bar{t}]$. Hence, every point $(x,t)$ with $(x_3,t) \in Q$ lies at the center of an $\varepsilon_0$-neck. 

We next consider the evolution equation of $rr_z$. We compute 
\[(rr_z)_t=(rr_z)_{zz}-\frac{(2+3r_z^2)r_zr_{zz} + rr_z^2r_{zzz}}{1+r_z^2} - \frac{2rr_zr_{zz}^2}{(1+r_z^2)^2}.\]
Using Lemma \ref{bound.for.rr_z} and Lemma \ref{interior.estimate.for.r}, we obtain  
\[\sup_Q |rr_z| \leq C\]
and
\[\sup_Q |(rr_z)_{t} - (rr_z)_{zz}| \leq C \bar{r}^{-3}.\] 
Standard interior estimates for the linear heat equation give 
\[|(rr_z)_z| \leq C \bar{r}^{-\frac{3}{2}} \sup_Q |rr_z| + C \bar{r}^{\frac{3}{2}} \sup_Q |(rr_z)_t - (rr_z)_{zz}| \leq C \bar{r}^{-\frac{3}{2}}\] 
at $(\bar{z},\bar{t})$. Thus, $|r_{zz}| \leq C \bar{r}^{-\frac{5}{2}}$ at $(\bar{z},\bar{t})$. This completes the proof of Proposition \ref{estimate.for.r_zz}. \\

For each $z<0$, we define a real number $\mathcal{T}(z)$ by
\begin{align*}
&r(z,t)>0 \quad \text{for}\quad t < \mathcal{T}(z), && \lim_{t \to \mathcal{T}(z)} r(z,t)=0.
\end{align*}
The following result allows us to estimate $r(z,t)$ in terms of $\mathcal{T}(z)-t$.

\begin{corollary}
\label{cor}
We have 
\[2[\mathcal{T}(z)-t] \leq r(z,t)^2 \leq 2[\mathcal{T}(z)-t] + 8C_2[\mathcal{T}(z)-t]^{\frac{1}{4}}+C_1^2\]
if $z<0$ and $r(z,t)$ is sufficiently large.
\end{corollary}

\textbf{Proof.} We again fix a point $(\bar{z},\bar{t})$. Since $(r^2+2t)_t = \frac{2rr_{zz}}{1+r_z^2} < 0$, we have 
\[r(\bar{z},\bar{t})^2 \geq 2[\mathcal{T}(\bar{z})-\bar{t}].\] 
Moreover, Proposition \ref{estimate.for.r_zz} implies that $(r^2+2t)_t = \frac{2rr_{zz}}{1+r_z^2} \geq -2C_2 r^{-\frac{3}{2}}$ whenever $r \geq C_1$. Let $\tilde{t} \leq \mathcal{T}(\bar{z})$ denote the time when $r(\bar{z},\tilde{t})=C_1$. Since $r(\bar{z},t)$ is a decreasing function of $t$, $r(\bar{z},t) \leq C_1$ for $t \leq \tilde{t}$. Therefore,
\begin{align*}
r(\bar{z},\bar{t})^2
&= C_1^2+2(\tilde{t}-\bar{t}) - \int_{\bar{t}}^{\tilde{t}} (r(\bar{z},t)^2+2t)_t \, dt \\
&\leq C_1^2+2(\tilde{t}-\bar{t})+2C_2\int_{\bar{t}}^{\tilde{t}} r(\bar{z},t)^{-\frac{3}{2}} \, dt \\ 
&\leq C_1^2+2(\tilde{t}-\bar{t})+2C_2\int_{\bar{t}}^{\tilde{t}} [\mathcal{T}(\bar{z})-t]^{-\frac{3}{4}} \, dt \\
&\leq C_1^2+2(\tilde{t}-\bar{t})-8C_2[\mathcal{T}(\bar{z})-\tilde{t}]^{\frac{1}{4}}+8C_2[\mathcal{T}(\bar{z})-\bar{t}]^{\frac{1}{4}} \\ 
&\leq C_1^2 + 2 \, [\mathcal{T}(\bar{z})-\bar{t}] + 8C_2[\mathcal{T}(\bar{z})-\bar{t}]^{\frac{1}{4}},
\end{align*} 
as claimed. \\

Lemma \ref{bound.for.rr_z} gives a sharp upper bound for the limit of $rr_z$. More precisely, $\limsup_{z \to \infty} r(z,t)r_z(z,t) \leq \mathcal{H}^{-1}$ for each $t$. We next establish a lower bound for $\liminf_{z \to \infty} r(z,t)r_z(z,t)$. To derive this estimate, we need a lower bound for $r(0,t) r_z(0,t)$.

\begin{lemma}
\label{lower.bound.for.rr_z.at.z=0}
Let $\delta>0$ be arbitrary. Then 
\[r(0,t)r_z(0,t) \geq \mathcal{H}^{-1}-\delta\]
whenever $-t$ is sufficiently large.
\end{lemma}

\textbf{Proof.} 
In the following, we assume that $-t$ is sufficiently large so that $R=r(0,t) \geq C_1$. Consequently, every point $(x,t)$ with $x_3=0$ lies at the center of an $\varepsilon_0$-neck. This implies $|r(z,t)-R| \leq \varepsilon_0 R$ for $|z| \leq 2R$. Recall that $rr_z \leq (1+2\varepsilon_0)\mathcal{H}^{-1}$ by Lemma \ref{bound.for.rr_z}, and $|(rr_{z})_z| = |rr_{zz}+r_z^2| \leq C_3R^{-\frac{3}{2}}$ for some constant $C_3$ by Proposition \ref{estimate.for.r_zz}. Hence, if we choose $-t$ sufficiently large so that $R^{\frac{1}{2}}\geq 4C_3\delta^{-1}$, then we obtain
\[|r(z,t)r_z(z,t)-r(0,t)r_z(0,\bar{t})| \leq 2C_3R^{-\frac{1}{2}} \leq \frac{\delta}{2}\]
for all $z \in [-2R,2R]$.

It follows from Corollary \ref{cor} that  
\[r(-R,t)^2 \geq 2[\mathcal{T}(-R)-t],\] 
\[r(-2R,t)^2 \geq 2[\mathcal{T}(-2R)-t],\] 
and 
\begin{align*} 
r(-2R,t)^2 
&\leq 2[\mathcal{T}(-2R)-t] + 8C_2 [\mathcal{T}(-2R)-t]^{\frac{1}{4}}+C_1^2 \\ 
&\leq 2[\mathcal{T}(-2R)-t] + 8C_2 r(-2R,t)^{\frac{1}{2}}+C_1^2 \\  
&\leq 2[\mathcal{T}(-2R)-t] + 8C_2 R^{\frac{1}{2}}+C_1^2. 
\end{align*}
This implies 
\[r(-R,t)^2-r(-2R,t)^2 \geq 2[\mathcal{T}(-R)-\mathcal{T}(-2R)]-8C_2R^{\frac{1}{2}}-C_1^2.\] 
Moreover, if $R$ is sufficiently large, then 
\[\mathcal{T}(-R) - \mathcal{T}(-2R) \geq \Big ( \mathcal{H}^{-1}-\frac{\delta}{2} \Big ) \, R.\] 
This gives 
\[r(-R,t)^2-r(-2R,t)^2 \geq 2 \, \Big ( \mathcal{H}^{-1}-\frac{\delta}{2} \Big ) \, R\] 
if $-t$ is sufficiently large. Hence, if $-t$ is sufficiently large, then 
\[\sup_{z \in [-2R,R]} r(z,t) r_z(z,t) \geq \mathcal{H}^{-1}-\frac{\delta}{2}.\] 
Putting these facts together, we conclude that  
\[r(0,t)r_z(0,t) \geq \mathcal{H}^{-1}-\delta\] 
whenever $-t$ is sufficiently large. \\

We next recall a solution $\psi(z,t)$ to the heat equation satisfying Dirichlet boundary condition on the half line. 

\begin{proposition}
\label{1d.heat.equation}
We define a smooth function $\psi: (0,\infty) \times (0,\infty) \to \mathbb{R}$ by
\begin{align*}
\psi(z,t)=\frac{1}{\sqrt{4\pi t}}\int^{\infty}_0 (e^{-\frac{(z-y)^2}{4t}}-e^{-\frac{(z+y)^2}{4t}}) dy.
\end{align*}
Then $\psi$ is a solution to the heat equation $\psi_t=\psi_{zz}$. Moreover, for each $z>0$ and $t>0$  we have $\psi_{zz}(z,t) <0$ and
\begin{align*}
& \lim_{z \to 0} \psi(z,t) = 0, && \lim_{z \to \infty}\psi(z,t) = 1,  && \lim_{t \to 0}\psi(z,t)=1, && \lim_{t \to \infty}\psi(z,t)=0.
\end{align*}
\end{proposition}

\textbf{Proof.} We only need to show $\psi_{zz}<0$. Direct computations yield
\begin{align*}
\psi_{zz} 
&= \frac{1}{\sqrt{4\pi t}} \bigg [ -\int^{\infty}_0 \Big ( \frac{1}{2t}-\frac{(z-y)^2}{4t^2} \Big ) e^{-\frac{(z-y)^2}{4t}} \, dy + \int_0^\infty \Big ( \frac{1}{2t}-\frac{(z+y)^2}{4t^2} \Big ) e^{-\frac{(z+y)^2}{4t}} dy \bigg ] \\
&= \frac{1}{\sqrt{8\pi t^2}}\bigg[ -\int^{\infty}_{-\frac{z}{\sqrt{2t}}}(1-\xi^2)e^{-\frac{\xi^2}{2}}d\xi+\int^{\infty}_{\frac{z}{\sqrt{2t}}}(1-\xi^2)e^{-\frac{\xi^2}{2}}d\xi\bigg]\\
&= -\frac{1}{\sqrt{8\pi t^2}} \int^{\frac{z}{\sqrt{2t}}}_{-\frac{z}{\sqrt{2t}}}(1-\xi^2)e^{-\frac{\xi^2}{2}}d\xi.
\end{align*}
Clearly, $\psi_{zz}<0$ for $0 < z \leq \sqrt{2t}$. Moreover, $\psi_{zzz} > 0$ for $z \geq \sqrt{2t}$, and 
\[\lim_{z \to \infty}\psi_{zz}(z,t)=-\frac{1}{\sqrt{8\pi t^2}} \int^\infty_{-\infty} (1-\xi^2) e^{-\frac{\xi^2}{2}} d\xi = 0.\] 
Therefore,  $\psi_{zz}<0$ also holds for $z \geq \sqrt{2t}$. This completes the proof of Proposition \ref{1d.heat.equation}. \\

\begin{proposition}
\label{liminf}
Given $\delta>0$, there exists a time $\bar{t} \in (-\infty,0]$ (depending on $\delta$) such that 
\begin{align*}
r(z,t)r_z(z,t) \geq \mathcal{H}^{-1}-2\delta,
\end{align*}
holds for all $z \geq 0$ and $t \leq \bar{t}$.
\end{proposition}

\textbf{Proof.} By Proposition \ref{estimate.for.r_zz}, we have $1+rr_{zz} \geq 0$ for $r \geq C_1+C_2$. This implies 
\[(rr_z)_t = \frac{(rr_z)_{zz}}{1+r_z^2}-\frac{2r_zr_{zz}(1+r_z^2+rr_{zz})}{(1+r_z^2)^2} \geq \frac{(rr_z)_{zz}}{1+r_z^2}\] 
for $r \geq C_1+C_2$. By Lemma \ref{lower.bound.for.rr_z.at.z=0}, we can choose $\bar{t}$ large enough so that $r(0,t)r_z(0,t) \geq \mathcal{H}^{-1}-\delta$ for $t \leq \bar{t}$. Moreover, by a suitable choice of $\bar{t}$ we can arrange that $r(z,t) \geq C_1+C_2$ for all $z \geq 0$ and $t \leq \bar{t}$. For each $s < \bar{t}$, we define a barrier function $\psi^{\delta,s}(z,t)$ by
\[\psi^{\delta,s}(z,t) = \mathcal{H}^{-1}-2\delta -\mathcal{H}^{-1} \, \psi(2z,t-s)\] 
for $t \in (s,\bar{t}]$. We claim that $rr_z > \psi^{\delta,s}$ for all $z \geq 0$ and all $t \in (s,\bar{t}]$.  

By our choice of $\bar{t}$, $r(0,t) r_z(0,t) \geq \mathcal{H}^{-1}-\delta > \limsup_{z \to 0} \psi^{\delta,s}(z,t)$ for each $t \in (s,\bar{t}]$. Moreover, $\liminf_{z \to \infty} r(z,t) r_z(z,t) \geq 0 > \limsup_{z \to \infty} \psi^{\delta,s} (z,t)$ for each $t \in (s,\bar{t}]$. Finally, Proposition \ref{1d.heat.equation} gives $r(z,s) r_z(z,s) \geq 0 > \limsup_{t \to s} \psi^{\delta,s} (z,t)$ for each $z>0$. 

Thus, if the inequality $rr_z > \psi^{\delta,s}$ fails, there exists some point $(z_0,t_0) \in (0,\infty) \times (s,\bar{t}]$ such that $r(z_0,t_0)r_z(z_0,t_0)=\psi^{\delta,s}(z_0,t_0)$ and $r(z,t)r_z(z,t) \geq \psi^{\delta,s}(z,t)$ for all $z \geq 0$ and all $t \in (s,t_0]$. Then, at the point $(z_0,t_0)$ we have
\[\frac{(\psi^{\delta,s})_{zz}}{1+r_z^2} \leq \frac{(rr_z)_{zz}}{1+r_z^2}\leq (rr_z)_t \leq (\psi^{\delta,s})_t = \frac{1}{4} \, (\psi^{\delta,s})_{zz}.\]
This contradicts the fact that $r_z \leq \varepsilon_0$ and $(\psi^{\delta,s})_{zz} >0$. 

Thus, we conclude that $rr_z >\psi^{\delta,s}$ for all $z \geq 0$ and all $t \in (s,\bar{t}]$. Sending $s \to -\infty$, we obtain $rr_z \geq \mathcal{H}^{-1}-2\delta$ for all $z \geq 0$ and all $t \leq \bar{t}$. \\ 

\begin{corollary}
\label{f.defined.everywhere}
We can find a time $T \in (-\infty,0]$ such that $r(z,t)^2 \geq \mathcal{H}^{-1} \, z$ for all $z \geq 0$ and $t \leq T$. In particular, if $t \leq T$, then the function $f(r,t)$ is defined for all $r \in [0,\infty)$.
\end{corollary}

\textbf{Proof.} By Proposition \ref{liminf}, we can find a time $T \in (-\infty,0]$ such that $r(z,t)r_z(z,t) \geq \frac{1}{2} \, \mathcal{H}^{-1}$ for all $z \geq 0$ and all $t \leq T$. From this, the assertion follows easily. \\

After these preparations, we now compute the limit $\lim_{z \to \infty} r(z,t)r_z(z,t)$. 

\begin{proposition}
\label{limit.of.rr_z}
For each $t \leq T$, we have $\lim_{z \to \infty} r(z,t)r_z(z,t) = \mathcal{H}^{-1}$.
\end{proposition}

\textbf{Proof.} Lemma \ref{bound.for.rr_z} gives $\limsup_{z\to \infty} r(z,t) r_z(z,t) \leq \mathcal{H}^{-1}$. So, it is enough to show that $\liminf_{z \to \infty} r(z,t) r_z(z,t) \geq \mathcal{H}^{-1}$ for each $t \leq T$. Given any $\delta>0$, Proposition \ref{liminf} implies that we can find a number $\bar{t} \leq T$ such that $\liminf_{z \to \infty} r(z,\bar{t}) r_z(z,\bar{t}) \geq \mathcal{H}^{-1} - 2\delta$. Moreover, Lemma \ref{interior.estimate.for.r} guarantees that 
\[|(rr_z)_t| = |r_tr_z+rr_{zt}| = \Big | \frac{r_zr_{zz}}{1+r_z^2}+\frac{rr_{zzz}}{1+r_z^2}-\frac{2rr_zr_{zz}^2}{(1+r_z^2)^2} \Big | \leq \frac{4C_0^3}{r^2}\] 
for $r \geq C_1$. Using Corollary \ref{f.defined.everywhere}, we obtain 
\[\liminf_{z \to \infty} r(z,t) r_z(z,t) = \liminf_{z \to \infty} r(z,\bar{t}) r_z(z,\bar{t}) \geq \mathcal{H}^{-1} - 2\delta\] 
for each $t \leq T$. Since $\delta>0$ is arbitrary, we conclude that $\liminf_{z \to \infty} r(z,t) r_z(z,t) \geq \mathcal{H}^{-1}$ for each $t \leq T$. This completes the proof of Proposition \ref{limit.of.rr_z}. \\

\begin{theorem}
For each $t \leq T$, $M_t$ is a translating soliton.
\end{theorem}

\textbf{Proof.} Since $rr_z =\frac{r}{f_r}$, Proposition \ref{limit.of.rr_z} implies
\[\lim_{r \to \infty} \frac{f_r(r,t)}{r}=\mathcal{H}\]
for each $t \leq T$. Using the evolution equation for $f(r,t)$, we obtain 
\[\lim_{r \to \infty} f_t(r,t)=\lim_{r\to \infty} \frac{f_r(r,t)}{r}=\mathcal{H}\] 
for each $t \leq T$. Using Proposition \ref{monotonicity.in.r}, we conclude that $f_t(r,t) \leq \mathcal{H}$ for all $r \geq 0$ and all $t \leq T$. Therefore, Proposition \ref{lower.bound.for.f_t} gives $f_t(r,t)=\mathcal{H}$ for all $r \geq 0$ and all $t \leq T$. Consequently, $M_t$ is a translating solition for each $t \leq T$. \\

Once we know that $M_t$ is a translating soliton for $-t$ sufficiently large, it follows from standard arguments that $M_t$ is a translating soliton for all $t$:

\begin{proposition}
Suppose that $M_{\bar{t}}$ is a translating soliton for some $\bar{t} \in (-\infty,0]$. Then $M_t$ is a translating soliton for all $t \in [\bar{t},0]$.
\end{proposition}

\textbf{Proof.} 
By Proposition \ref{H_max.finite} and Corollary \ref{properties.of.H_max}, the flow $M_t$, $t \in (-\infty,0]$, has bounded curvature. Without loss of generality, we may assume that $\sup_{M_t} |A|^2 \leq 2$ for each $t \in (-\infty,0]$. If $\omega$ is a fixed vector in $\mathbb{R}^3$, then 
\[\frac{\partial}{\partial t} (H-\langle \omega,\nu \rangle) = \Delta (H-\langle \omega,\nu \rangle) + |A|^2 (H-\langle \omega,\nu \rangle).\] 
Moreover, since $|A|^2 \leq 2$, the function $\rho(x,t) := e^{8t} (|x|^2+1)$ satisfies 
\[\frac{\partial}{\partial t} \rho > \Delta \rho + |A|^2 \rho\] 
for $t \in (-\infty,0]$. By the maximum principle, the quantity $\sup_{M_t} \frac{|H - \langle \omega,\nu \rangle|}{\rho}$ is monotone decreasing for $t \in (-\infty,0]$. In particular, if $H = \langle \omega,\nu \rangle$ at each point on $M_{\bar{t}}$, then $H = \langle \omega,\nu \rangle$ on $M_t$ for all $t \in [\bar{t},0]$. \\

\end{document}